\documentclass[onefignum,onetabnum]{siamart190516}

\usepackage{lipsum}
\usepackage{amsfonts}
\usepackage[square, numbers]{natbib}
\usepackage{graphicx}
\usepackage{epstopdf}
\usepackage{algorithmic}
\usepackage{float}
\usepackage{graphicx}

\usepackage{blindtext}
\usepackage{bm, bbm}
\usepackage{standalone}
\usepackage[misc]{ifsym}
\usepackage{amsmath, amssymb}
\usepackage{mathtools}
\usepackage{pgfplots}
\usepackage{stmaryrd}
\usepackage{comment}
\usepackage{xfrac}
\usepackage{url}
\usepackage{amsmath}
\usepackage{algorithm}
\usepackage{algorithmic}
\usepackage{subcaption}
\usepackage{booktabs}
\usepackage{multirow}
\usepackage{placeins}
\usepackage{accents}
\usepackage{graphics}

\usepackage{mathrsfs}
\usepackage{tikz}
\ifpdf
  \DeclareGraphicsExtensions{.eps,.pdf,.png,.jpg}
\else
  \DeclareGraphicsExtensions{.eps}
\fi

\usepackage{mdframed}
\mdfsetup{skipabove=22pt,skipbelow=22pt}
\newmdenv[
  topline=true,
  bottomline=true,
  rightline=true,
  leftline=true,
  innertopmargin=11pt,
  linewidth=.75pt
]{exbox}
\newenvironment{box_example}{\begin{exbox}\begin{example}}{\end{example}\end{exbox}}

\newsiamremark{remark}{Remark}
\newsiamremark{hypothesis}{Hypothesis}
\crefname{hypothesis}{Hypothesis}{Hypotheses}
\newsiamthm{claim}{Claim}
\newsiamthm{example}{Example}
\newsiamthm{assumption}{Assumption}

\headers{A Unified Approach to Mixed-Integer Optimization}{D. Bertsimas, R. Cory-Wright, and J. Pauphilet}

\title{A Unified Approach to Mixed-Integer Optimization Problems With Logical Constraints}

\author{Dimitris Bertsimas\thanks{Sloan School of Management, Massachusetts Institute of Technology, Cambridge, MA, USA (\email{dbertsim@mit.edu}).}
\and Ryan Cory-Wright\thanks{Operations Research Center, Massachusetts Institute of Technology, Cambridge, MA, USA (\email{ryancw@mit.edu}).}
\and Jean Pauphilet\thanks{London Business School, London, UK
  (\email{jpauphilet@london.edu}).}}

\usepackage{amsopn}

\makeatletter
\newcommand*{\addFileDependency}[1]{% argument=file name and extension
  \typeout{(#1)}% latexmk will find this if $recorder=0 (however, in that case, it will ignore #1 if it is a .aux or .pdf file etc and it exists! if it doesn't exist, it will appear in the list of dependents regardless)
  \@addtofilelist{#1}% if you want it to appear in \listfiles, not really necessary and latexmk doesn't use this
  \IfFileExists{#1}{}{\typeout{No file #1.}}% latexmk will find this message if #1 doesn't exist (yet)
}
\makeatother

\usepackage{color}

\ifpdf
\hypersetup{
  pdftitle={A Unified Approach to Mixed-Integer Optimization Problems With Logical Constraints},
  pdfauthor={D. Bertsimas, R. Cory-Wright, and J. Pauphilet}
}
\fi

\usepackage{xcolor}
\hypersetup{
    colorlinks,
    linkcolor={red!50!black},
    citecolor={blue!50!black},
    urlcolor={blue!80!black}
}
\begin{document}

\maketitle

% REQUIRED
\begin{abstract}
We propose a unified framework to address a family of classical mixed-integer optimization problems with {logically constrained} decision variables, including network design, facility location, unit commitment, sparse portfolio selection, binary quadratic optimization{, sparse principal component analysis,} and sparse learning problems. These problems exhibit logical relationships between continuous and discrete variables, which are usually reformulated linearly using a big-$M$ formulation. In this work, we challenge this longstanding modeling practice and express the logical constraints in a non-linear way. By imposing a regularization condition, we reformulate these problems as convex binary optimization problems, which are solvable using an outer-approximation procedure.
In numerical experiments, we establish that a general-purpose numerical strategy, which combines cutting-plane, first-order, and local search methods, solves these problems faster and at a larger scale
than state-of-the-art mixed-integer linear or second-order cone methods.
Our approach successfully solves network design problems with $100$s of nodes and provides solutions up to $40\%$ better than the state-of-the-art; sparse portfolio selection problems with up to $3,200$ securities compared with $400$ securities for previous attempts; and sparse regression problems with up to $100,000$ covariates.
\end{abstract}

% REQUIRED

% REQUIRED

\section{Introduction}\label{sec:intro}
Many important problems from the Operations Research literature exhibit a logical relationship between continuous variables $x$ and binary variables $z$ of the form ``$x=0$ if $z=0$''. Among others, start-up costs in machine scheduling problems, financial transaction costs, cardinality constraints and fixed costs in facility location problems exhibit this relationship. Since the work of \cite{glover1975improved}, this relationship is usually enforced through a ``big-$M$'' constraint of the form ${-M z \leq x \leq M z}$ for a sufficiently large constant $M >0$. Glover's work has been so influential that big-$M$ constraints are now considered as intrinsic components of the initial problem formulations themselves, to the extent that textbooks in the field introduce facility location, network design or sparse portfolio problems with big-$M$ constraints \emph{by default}, although they are actually \emph{reformulations} of logical constraints.

In this work, we adopt a different perspective on the big-$M$ paradigm, viewing it as a regularization term, rather than a modeling trick. Under this lens, we show that regularization drives the computational tractability of problems with logical constraints, explore alternatives to the big-$M$ paradigm and propose an efficient algorithmic strategy which solves a broad class of problems with logical constraints.

\subsection{Problem Formulation and Main Contributions}
We consider optimization problems which unfold over two stages. In the first stage, a decision-maker activates binary variables, while satisfying resource budget constraints and incurring activation costs. Subsequently, in the second stage, the decision-maker optimizes over the continuous variables. Formally, we consider the problem
\begin{equation}
\begin{aligned}\label{eqn:original_minlp}
    \min_{\bm{z} \in \mathcal{Z}, {\bm{x} \in  \mathbb{R}^n}} \quad & \bm{c}^\top \bm{z} + g(\bm{x}) + \Omega(\bm{x}) \quad
    \text{s.t.} \quad x_i =0\ \text{if} \  z_i=0 \quad \forall i \in [n],
\end{aligned}
\end{equation}
where $\mathcal{Z} \subseteq \{0,1\}^n$, $\bm{c} \in \mathbb{R}^n$ is a cost vector, $g(\cdot)$ is a generic convex function {which possibly models convex constraints $\bm{x} \in \mathcal{X}$ for a convex set $\mathcal{X} \subseteq \mathbb{R}^n$ implicitly—by requiring that $g(\bm{x})=+\infty$ if $\bm{x} \notin \mathcal{X}$}, and $\Omega(\cdot)$ is a convex regularization function; we formally state its structure in Assumption \ref{assump:sep}.

In this paper, we provide three main contributions: First, we reformulate the logical constraint ``$x_i=0$ if $z_i=0$'' in a non-linear way, by substituting $z_i x_i$ for $x_i$  in Problem \eqref{eqn:original_minlp}. Second, we leverage the
%decomposition of the objective function into a convex function $g(\bm{x})$ plus a
regularization term $\Omega(\bm{x})$ to derive a tractable reformulation of \eqref{eqn:original_minlp}. Finally, by invoking strong duality, we reformulate \eqref{eqn:original_minlp} as a mixed-integer saddle-point problem, which is solvable via outer approximation.

Observe that the structure of Problem \eqref{eqn:original_minlp} is quite general, as the feasible set $\mathcal{Z}$ can capture known lower and upper bounds on $\bm{z}$, relationships between different $z_i$'s, or a cardinality constraint $ \bm{e}^\top \bm{z} \leq k$. Moreover, constraints of the form $\bm{x} \in \mathcal{X}$, for some convex set $\mathcal{X}$, can be encoded within the domain of $g$, by defining $g(\bm{x})=+ \infty$ if $\bm{x} \notin \mathcal{X}$. As a result, Problem \eqref{eqn:original_minlp} encompasses a large number of problems from the Operations Research literature, such as the network design problem described in Example \ref{ex:nd}. These problems are typically studied separately. %, resulting in the development of problem-specific strategies. H
However, the techniques developed for each problem are actually different facets of a single unified story, and, as we demonstrate in this paper, can be applied to a much more general class of problems than is often appreciated.

\begin{box_example} \label{ex:nd}
\textit{Network design is an important example of problems of the form \eqref{eqn:original_minlp}. Given a set of $m$ nodes, the network design problem consists of constructing edges to minimize the construction plus flow transportation cost. Let $E$ denote the set of all potential edges and let $n = |E|$. Then, the network design problem is given by:
\begin{equation} \label{eqn:nd} %\tag{ND}
\begin{aligned}
    \min_{\bm{z} \in \mathcal{Z}, {\bm{x} \in \mathbb{R}^n_{+}}} \quad \bm{c}^\top \bm{z} + \tfrac{1}{2} \bm{x}^\top \bm{Q} \bm{x} + \bm{d}^\top \bm{x} \quad
 \mbox{ s.t. } & \ \bm{A} \bm{x} = \bm{b}, \\
    & \ x_e = 0 \mbox{ if } z_e = 0 \quad \forall e \in E,
\end{aligned}
\end{equation}
where {$\mathcal{Z} = \{0,1\}^n$}, $\bm{A} \in \mathbb{R}^{m \times n}$ is the flow conservation matrix, $\bm{b} \in \mathbb{R}^{m}$ is the vector of external demands and $\bm{Q} \in \mathbb{R}^{n \times n}$, $\bm{d} \in \mathbb{R}^{n}$ define the quadratic and linear costs of flow circulation. We assume that $\bm{Q} \succeq \bm{0}$ is a positive semidefinite matrix. Inequalities of the form $\bm{\ell} \leq \bm{z} \leq \bm{u}$ can be incorporated within $\mathcal{Z}$ to account for existing/forbidden edges in the network. Problem \eqref{eqn:nd} is of the same form as Problem \eqref{eqn:original_minlp} with
\begin{align*}
    g(\bm{x}) + \Omega(\bm{x}) := \begin{dcases}
    \tfrac{1}{2} \bm{x}^\top \bm{Q} \bm{x} + \bm{d}^\top \bm{x}, &\mbox{ if } \bm{A} \bm{x} = \bm{b}, \bm{x} \geq \bm{0}, \\
    +\infty, & \mbox{ otherwise}.
    \end{dcases}
\end{align*}
We present a generalized model with edge capacities and  multiple commodities in Section \ref{sec:ex_nd}.}
\end{box_example}
\subsection{Background and Literature Review}\label{ssec:background}
Our work falls into two areas of the mixed-integer optimization literature which are often considered in isolation: $(a)$ modeling forcing constraints which encode whether continuous variables are active and can take non-zero values or are inactive and forced to $0$, and $(b)$ decomposition algorithms for mixed-integer optimization problems.

\subsubsection*{Formulations of forcing constraints} %First part on how to formulate the logical constraints to fall into solvable family of optimization problem
The most popular way to impose forcing constraints on continuous variables is to introduce auxiliary discrete variables which encode whether the continuous variables are active, and relate the discrete and continuous variables via
%enforce that discrete variables encode whether continuous variables are activated is
the big-$M$ approach of \cite{glover1975improved}. This approach was first applied to mixed-integer non-linear optimization (MINLO) in the context of sparse portfolio selection by \cite{bienstock1996computational}.  With the big-$M$ approach, the original MINLO admits bounded relaxations %problem becomes bounded,
%with additional mixed-integer linear constraints,
and can therefore be solved via branch-and-bound.
%In addition, since the constraints between the discrete and continuous variables are linear, they have a theoretically low impact on the computational tractability of the optimization problems.
Moreover, because the relationship between discrete and continuous variables is enforced via linear constraints, a big-$M$ reformulation has a theoretically low impact on the tractability of the MINLOs continuous relaxations.
However, in practice, high values of $M$ lead to numerical instability and provide low-quality bounds \citep[see][Section 5]{beaumont1990algorithm}.

This observation led \cite{frangioni2006perspective} to propose a class of cutting-planes for MINLO problems with indicator variables, called perspective cuts, which often provide a tighter reformulation of the logical constraints. Their approach was subsequently extended by \cite{akturk2009strong}, who{, building upon the work of \cite[pp. 88, item 5]{ben2001lectures}}, proved that MINLO problems with indicator variables can often be reformulated as mixed-integer second-order cone problems (see \cite{gunluk2012perspective} for a survey). More recently, a third approach for coupling the discrete and the continuous in MINLO was proposed independently for sparse regression by \cite{pilanci2015sparse} and \cite{bertsimas2017sparse}: augmenting the objective with a strongly convex term of the form $\Vert\bm{x}\Vert_2^2$, called a ridge regularizer.

In the present paper, we synthesize the aforementioned and seemingly unrelated three lines of research under the unifying lens of regularization. Notably, our framework includes big-$M$ and ridge regularization as special cases, and provides an elementary derivation of perspective cuts.

\subsubsection*{Numerical algorithms for mixed-integer optimization} %Second part on how to solve large scale MIO

A variety of ``classical'' general-purpose decomposition algorithms have been proposed for general MINLOs. The first such decomposition method is known as
{Generalized Benders Decomposition, and was proposed by \cite{geoffrion1972generalized} %, who proved its finite termination, building upon the work of
as an extension of \cite{benders1962partitioning}. A similar method, known as}
outer-approximation was proposed by
%Benders decomposition, and was proposed for linear optimization problems by \cite{benders1962partitioning}.
\cite{duran1986outer}, who proved its finite termination. The outer-approximation method was subsequently generalized to account for non-linear integral variables by \cite{fletcher1994solving}. These techniques decompose MINLOs into a discrete master problem and a sequence of continuous separation problems, which are iteratively solved to generate valid cuts for the master problem.% on-the-fly.

Though slow in their original implementation, decomposition schemes have benefited from recent improvements in mixed-integer linear solvers in the past decades, beginning with the
% Unfortunately however, ``classical'' decomposition schemes such as outer-approximation converge notoriously slowly in practise.
% The slow convergence of ``classical'' outer-approximation has motivated the development of ``modern'' decomposition schemes, beginning with the
%a tremendous number of  over the past decades
branch-and-cut approaches of \citep{padberg1991branchcut, quesada1992lp}, which embed the cut generation process within a single branch-and-bound tree, rather than building a branch-and-bound tree before generating each cut. We refer to \cite{fischetti2016benders, fischetti2016redesigning} for recent successful implementations of ``modern'' decomposition schemes. From a high-level perspective, these recent successes require three key ingredients: First, a fast cut generation strategy. Second, as advocated by \cite{fischetti2016benders}, a rich cut generation process at the root node. Finally, a cut selection rule for degenerate cases where multiple valid inequalities exist (e.g., the Pareto optimality criteria of \cite{magnanti1981accelerating}).

In this paper, we connect the regularization used to reformulate logical constraints with %constitutes a modeling choice, and connect this choice to
the aforementioned key ingredients for modern decomposition schemes. %This approach allows us to
%connect the modeling choice to be made in reformulating the logical constraints with such numerical considerations, and
Hence, instead of considering a MINLO formulation as a given and %we also argue that rather than making a modeling choice to fix a MINLO formulation and
subsequently
%shift our focus from
attempt to solve it at scale,
our approach view big-$M$ constraints as one of many alternatives. We argue that regularization is a modeling choice that impacts the tractability of the formulation and should be made accordingly.

\subsection{Structure}
We propose a unifying framework to address mixed-integer optimization problems, and jointly discuss modeling choice and numerical algorithms.

In Section \ref{sec:framework}, we identify a general class of mixed-integer optimization problems, which encompasses sparse regression, sparse portfolio selection, sparse principal component analysis, unit commitment, facility location, network design and binary quadratic optimization as special cases. For this class of problems, we discuss how imposing either big-$M$ or ridge regularization accounts for non-linear relationships between continuous and binary variables in a tractable fashion. We also establish that regularization controls the convexity and smoothness of Problem \eqref{eqn:original_minlp}'s objective function.

In Section \ref{sec:algo}, we propose a conjunction of general-purpose numerical algorithms to solve Problem \eqref{eqn:original_minlp}. The backbone of our approach is an outer approximation framework, enhanced with first-order methods to solve the Boolean relaxations and obtain improved lower bounds, certifiably near-optimal warm-starts via randomized rounding, and a discrete local search procedure. We also connect our approach to the perspective cut approach \citep{frangioni2006perspective} from a theoretical and implementation standpoint.

Finally, in Section \ref{sec:numresults}, we demonstrate empirically that algorithms derived from our framework can outperform state-of-the-art solvers. On network design problems with $100$s of nodes and binary quadratic optimization problems with $100$s of variables, we improve the objective value of the returned solution by $5$ to $40\%$ and $5$ to $85\%$ respectively, and our edge increases as the problem size increases. On empirical risk minimization problems, our method with ridge regularization is able to accurately select features among $100,000$s (resp. $10,000$s) of covariates for regression (resp. classification) problems, with higher accuracy than both Lasso and non-convex penalties from the statistics literature. For sparse portfolio selection, we solve to provable optimality problems one order of magnitude larger than previous attempts. We then analyze the benefits of the different ingredients in our numerical recipe on facility location problems, and discuss the relative merits of different regularization approaches on unit commitment instances.

\subsection*{Notation}
We use nonbold face characters to denote scalars {and components of matrices}, lowercase bold faced characters such as $\bm{x}$ to denote vectors, uppercase bold faced characters such as $\bm{X}$ to denote matrices, and calligraphic characters such as $\mathcal{X}$ to denote sets. We let $\mathbf{e}$ denote a vector of all $1$'s, and $\bm{0}$ denote a vector of all $0$'s, with dimension implied by the context. %We let $\succeq$ denote the L{\"o}wner partial order on the space of real symmetric matrices.  %If $x$ is a real number then $\lfloor x \rceil$ denotes the closest integer to $x$.
If $\bm{x}$ is a $n$-dimensional vector then $\mathrm{Diag}(\bm{x})$ denotes the $n \times n$ diagonal matrix whose diagonal entries are given by $\bm{x}$. If $f(\bm{x})$ is a convex function then its perspective function $\varphi(\bm{x}, t)$, defined as $\varphi(\bm{x}, t) = tf(\bm{x}/t)$ if $t>0$, $\varphi(\mathbf{0}, 0) =0$, and $\infty$ elsewhere, is also convex \citep[Chapter 3.2.6.]{boyd2004convex}. Finally, we let $\mathbb{R}_+^n$ denote the $n$-dimensional nonnegative orthant.
%Finally, we let $x_{[k]}$ denote the $k$th largest entry of the vector $\bm{x}$. %We denote the number of non-zero elements in a vector by the support norm $\ell_0$.

\section{Framework and Examples}\label{sec:framework}
In this section, we present the family of problems to which our analysis applies, discuss the role played by regularization, and provide some examples from the Operations Research literature.

\subsection{Examples} \label{sec:examples}
Problem \eqref{eqn:original_minlp} has a two-stage structure which comprises first ``turning on'' some indicator variables $\bm{z}$, and second solving a continuous optimization problem over the active components of $\bm{x}$. Precisely, Problem \eqref{eqn:original_minlp} can be viewed as a discrete optimization problem:
\begin{align}\label{eqn:masterproblem}
    \min_{\bm{z} \in \mathcal{Z}} \quad & \bm{c}^\top \bm{z} + f(\bm{z}),
\end{align}
where the inner minimization problem
\begin{equation}
\begin{aligned}\label{eqn:subproblem}
  f(\bm{z}) :=\min_{\bm{x} \in \mathbb{R}^n} &\quad  g(\bm{x}) + \Omega(\bm{x}) \quad
   \text{s.t.} \quad x_i =0\ \text{if}\  z_i=0  \quad \forall i \in [n],
\end{aligned}
\end{equation}
yields {a} best choice of $\bm{x}$ given $\bm{z}$. As we illustrate in this section, a number of problems of practical interest exhibit this structure.\FloatBarrier
\begin{box_example} \textit{
For the network design example \eqref{eqn:nd}, we have
\begin{align*}
  f(\bm{z}) :=\min_{{\bm{x} \in \mathbb{R}^n_+} : \bm{A} \bm{x} = \bm{b}} &\quad \tfrac{1}{2} \bm{x}^\top \bm{Q} \bm{x} + \bm{d}^\top \bm{x} \quad \text{s.t.} \quad x_e =0\ \text{if}\  z_e=0  \quad \forall e \in E.
\end{align*}}
\end{box_example}
\subsubsection{Network Design} \label{sec:ex_nd}
Example \ref{ex:nd} illustrates that the single-commodity network design problem is a special case of Problem \eqref{eqn:original_minlp}. We now formulate the $k$-commodity network design problem with directed capacities {as minimizing over  $\mathcal{Z} = \{0,1\}^n$ the function}:
\begin{equation}\label{eqn:mcnd}
\begin{aligned}
    f(\bm{z}) :=\min_{{\bm{f}^j, \bm{x} \in \mathbb{R}^n_+}} \quad \tfrac{1}{2} \bm{x}^\top \bm{Q} \bm{x} + \bm{d}^\top \bm{x} \quad
 \mbox{ s.t. } & \ \bm{A} \bm{f}^j = \bm{b}^j\quad \forall j\in[k], \\
 & \ \bm{x} = \sum_{j=1}^{k} \bm{f}^j, \ \bm{x} \leq \bm{u}, \\
    & \ x_e = 0 \mbox{ if } z_e = 0 \quad \forall e \in E.
\end{aligned}
\end{equation}

\subsubsection{Sparse Empirical Risk Minimization}
Given a matrix of covariates $\bm{X} \in \mathbb{R}^{n \times p}$ and a response vector $\bm{y} \in \mathbb{R}^n$, the sparse empirical risk minimization problem seeks
a vector $\bm{w}$ which explains the response in a compelling manner, i.e., {minimizes over $\mathcal{Z}:=\{\bm{z} \in \{0, 1\}^p: \bm{e}^\top \bm{z} \leq k\}$ the function:}
\begin{equation} \label{eqn:serm}
\begin{aligned}
    f(\bm{z}):=\min_{\bm{w} \in \mathbb{R}^p} \ & \sum_{i=1}^n \ell\left(y_i, \bm{w}^\top \bm{x}_i\right)+\frac{1}{2\gamma} \| \bm{w} \|_2^2\quad \text{s.t.} \quad w_{j} =0\ \text{if}\  z_{j}=0 \quad \forall j \in [p],
\end{aligned}
\end{equation}
where $\ell$ is an appropriate convex loss function; we provide examples of suitable loss functions in Table \ref{tab:ermlossfunctions_fenchel}.
\begin{table}[h!]
\centering
\caption{Loss functions and Fenchel conjugates for ERM problems of interest.}
\begin{tabular}{llll} \toprule
Method & Loss function & Domain & Fenchel conjugate\\\midrule
OLS & $\frac{1}{2}(y-u)^2$ & $y \in \mathbb{R}$ & $\ell^\star(y, \alpha) = \frac{1}{2} \alpha^2 + \alpha y$\\
SVM & $\max(1 - yu, 0)$ & $y \in \{\pm 1\}$ & $\ell^\star(y, \alpha) =
\begin{cases}
\alpha y, \ \text{if} \ \alpha y \in [-1,0],\\
\infty, \ \text{otherwise}.
\end{cases}$ \\
\bottomrule
\end{tabular}
\label{tab:ermlossfunctions_fenchel}
\end{table}
\subsubsection{Sparse Portfolio Selection}
Given an expected marginal return vector $\bm{\mu} \in \mathbb{R}^n$, estimated covariance matrix ${\bm{\Sigma} \in \mathcal{S}^n_+}$, uncertainty budget parameter $\sigma >0$, cardinality budget parameter $k \in \{2, \ldots, n-1\}$, linear constraint matrix $\bm{A} \in \mathbb{R}^{n \times m}$, and right-hand-side bounds $\bm{l}, \bm{u} \in \mathbb{R}^m$, investors determine an optimal allocation of capital between assets by minimizing over $\mathcal{Z} = \left\{ \bm{z} \in \{0, 1\}^n: \ \bm{e}^\top \bm{z} \leq k \right\}$ the function
\begin{equation}  \label{eqn:sps} %\tag{SPS}
\begin{aligned}
   f(\bm{z}):= \ \min_{\bm{x} \in \mathbb{R}_+^n}\quad & \frac{\sigma}{2} \bm{x}^\top \bm{\Sigma} \bm{x} - \bm{\mu}^\top \bm{x} \\
    \text{s.t.} \quad &  \bm{l} \leq \bm{A}\bm{x} \leq \bm{u},\ \bm{e}^\top \bm{x}=1, \ x_{i} =0\ \text{if} \  z_{i}=0\quad \forall i \in [n].
\end{aligned}
\end{equation}

\subsubsection{Unit Commitment} \label{sec:ex_uc}
In the DC-load-flow unit commitment problem, each generation unit $i$ incurs a cost given by a quadratic cost function $f^i(x)=a_i x^2+b_i x+c_i$ for its power generation output $x \in [0, u_i]$.
Let $\mathcal{T}$ denote a finite set of time periods covering a time horizon (e.g., $24$ hours).  At each time period $t \in \mathcal{T}$, there is an estimated demand $d_{t}$. The objective is to generate sufficient power to satisfy demand at minimum cost, while respecting minimum time on/time off constraints. %(see Figure \ref{fig:UC_cost} for an example).
% \begin{figure}[h]
%     \centering
%     \includegraphics{perspectivefunction_uc_gr1.eps}
%     \caption{The generator cost function $\theta(x,z)=10 x^2+45x+1500z$, $x \in [0, 100z], z \in \{0, 1\}$ (blue) and its perspective relaxation (orange).}
%     \label{fig:UC_cost}
% \end{figure}

By introducing binary variables $z_{i,t}$, which denote whether generation unit $i$ is active in time period $t$, requiring that $\bm{z} \in \mathcal{Z}$, i.e., $\bm{z}$ obeys physical constraints such as minimum time on/off, the unit commitment problem admits the formulation:
\begin{align}
    \min_{\bm{z}} \quad & f(\bm{z}) + \sum_{t \in \mathcal{T}} \sum_{i=1}^n c_{i}z_{i,t}
   \quad \text{s.t.} \quad \bm{z} \in \mathcal{Z} \subseteq \{0,1\}^{n \times \vert \mathcal{T} \vert},
\end{align}
\begin{equation}\label{eqn:UC}
\begin{aligned}
   \text{where: }\quad f(\bm{z}):=\min_{\bm{x}} \quad & \sum_{t \in \mathcal{T}} {\left( \sum_{i=1}^n \tfrac{1}{2 }a_i x_{i,t}^2+b_i x_{i,t} \right) } \ \text{s.t.} \ \sum_{i=1}^n x_{i,t}\geq D_{t} \quad \forall t \in \mathcal{T}, \\
      & x_{i,t} \in [0, u_{i,t}] \quad \forall i \in [{n}], \forall t \in \mathcal{T},\\
    & x_{i,t} =0\ \text{if} \  z_{i,t}=0  \quad \forall i \in [{n}], \forall t \in \mathcal{T}.
\end{aligned}
\end{equation}

%Observe that a necessary and sufficient condition for $f(\bm{z})$ to be feasible for a given $\bm{z}$ is that $\sum_i z_{i,t} u_{i,t} \geq D_t, \forall t \in \mathcal{T}$. Therefore, we will impose this constraint in our first-stage problem a priori, in order to avoid the undesirable possibility of solving infeasible subproblems.
\subsubsection{Facility Location}
Given a set of $n$ facilities and $m$ customers, the facility location problem consists of constructing facilities $i \in [n]$ at cost $c_i$ to satisfy demand at minimal cost, i.e., {minimizing over $\mathcal{Z}=\{0, 1\}^n$ the function:}
\begin{equation} \label{eqn:flp} %\tag{FLP}
\begin{aligned}
    {f(\bm{z}):=}\ \min_{\bm{X}\in \mathbb{R}_+^{n \times m}} \: & \bm{c}^\top \bm{z} +  \sum_{j=1}^m \sum_{i=1}^n {c}_{ij} {x}_{ij} \
    \text{s.t.} \ \sum_{j=1}^m {x}_{ij} \leq {u}_{i} \quad \forall i \in [n], \\
    & \sum_{i=1}^n {x}_{ij} = d_j \quad \forall j \in [m], \ {x}_{ij} = 0 \ \text{ if } \  z_{i}=0 \quad \forall i \in [n], j \in [m].
\end{aligned}
\end{equation}
In this formulation, ${x}_{ij}$ corresponds to the quantity produced in facility $i$ and shipped to customer $j$ at a marginal cost of ${c}_{ij}$. Moreover, each facility $i$ has a maximum output capacity of ${u}_i$ and each customer $j$ has a demand of $d_j$. In the uncapacitated case where ${u}_{i} = \infty$, the inner minimization problems decouple into independent knapsack problems for each customer $j$.

{
\subsubsection{Sparse Principal Component Analysis (PCA)}
Given a {$p\times p$ positive semidefinite} covariance matrix $\bm{\Sigma}$, $\bm{\Sigma} \in S^p_+$ in short, the sparse PCA problem is to select a vector $\bm{z}$ which maximizes {over $\mathcal{Z} = \left\{ \bm{z} \in \{0, 1\}^p: \ \bm{e}^\top \bm{z} \leq k \right\}$ the function}
\begin{align}
    f(\bm{z})=\max_{\bm{x} \in \mathbb{R}^p} \quad \bm{x}^\top \bm{\Sigma}\bm{x} \ \text{s.t.} \ \Vert\bm{x}\Vert_2^2 = 1, x_i=0 \ \text{if} \ z_i=0 \quad \forall i \in [p].
\end{align}
This function is apparently non-concave in $\bm{z}$, because $f(\bm{z})$ is the optimal value of a non-convex quadratic optimization problem. Fortuitously however, this problem admits an exact mixed-integer semidefinite reformulation, namely
\begin{align}\label{sparsepcamisdp1}
    f(\bm{z})=\max_{\bm{X} \in S^p_+} \ \langle \bm{\Sigma}, \bm{X}\rangle \ \text{s.t.} \ \mathrm{tr}(\bm{X})=1, {x}_{i,j}=0 \ \text{if} \ z_i=0\ \text{or} \ z_j=0 \quad \forall i,j \in [p].
\end{align}
Indeed, for any fixed $\bm{z}$, Problem \eqref{sparsepcamisdp1} maximizes a linear function in $\bm{X}$ and therefore admits a rank-one optimal solution. Thus, we prove that sparse PCA admits an exact {mixed-integer semidefinite optimization} reformulation.

%By encoding the logical constraints with a big-$M$ penalty, Problem \eqref{sparsepcamisdp1}'s objective remains linear and some rank-one matrix $\bm{X}$ is optimal. Note, however, that this argument does not hold under ridge regularization. In this case, we can identify an optimal $\bm{z}$ under ridge regularization and subsequently re-optimize for $\bm{X}$. Since $$ %\frac{1}{2\gamma}\Vert \bm{X}\Vert_F^2 =
% \dfrac{1}{2\gamma}\sum_{i,j} X_{i,j}^2\leq \dfrac{1}{2\gamma} \left( \sum_i X_{i,i} \right)^2=\frac{1}{2\gamma},$$ where the inequality follows from the $2 \times 2$ minors in $\bm{X} \succeq \bm{0}$, the difference in objectives between the unregularized and regularized problems is at most $1/2\gamma$ and becomes negligible as $\gamma \rightarrow \infty$.
}

\subsubsection{Binary Quadratic Optimization}
Given a symmetric cost matrix $\bm{Q}$, the binary quadratic optimization problem consists of selecting a vector of binary variables $\bm{z}$ which {minimizes over $\mathcal{Z}=\{0, 1\}^n$ the function:}
\begin{equation}\label{eqn:bqp}
\begin{aligned}
    f(\bm{z})=\bm{z}^\top \bm{Q} \bm{z}.
\end{aligned}
\end{equation}

This formulation is non-convex and does not include continuous variables. However, introducing auxiliary continuous variables yields the equivalent formulation \citep{fortet1960applications} {of minimizing over $\mathcal{Z}=\{0, 1\}^n$ the function}:
\begin{align}
    {f(\bm{z}):=}\ \min_{\bm{Y} \in \mathbb{R}^{n \times n}_{+}} \quad  \langle \bm{Q}, \bm{Y} \rangle \quad \text{s.t.} \quad
    & {y}_{i,j} \leq 1 \quad & \forall i,j \in [n], \\
    & {y}_{i,j} \geq z_i+z_j-1 \quad & \forall i \in [n], \forall j \in [n]\backslash \{i\}, \nonumber\\
    & {y}_{i,i} \geq z_i \quad & \forall i \in [n], \nonumber\\
    & {y}_{i,j} = 0 \mbox{ if } z_i=0 \quad & \forall i,j \in [n],\nonumber\\
    & {y}_{i,j} = 0 \mbox{ if } z_j = 0 \quad & \forall i,j \in [n].\nonumber
\end{align}

% We remark that the triangle inequalities \citep[see][and references therein]{deza2009geometry} substantially improve the quality of binary quadratic relaxations and can easily be added via lazy callbacks within numerical solvers. We will make use of these inequalities in our numerical experiments and state them for completeness:
% \begin{align*}
%     Y_{i,j}+Y_{i,k}+Y_{j,k}-z_i-z_j-z_k \geq -1, \quad \forall i, j, k \in [n],\\
%     -Y_{i,j}-Y_{i,k}+Y_{j,k}+z_i \geq 0, \quad \forall i, j, k \in [n],\\
%     -Y_{i,j}+Y_{i,k}-Y_{j,k}+z_j \geq 0, \quad \forall i, j, k \in [n],\\
%     Y_{i,j}-Y_{i,k}-Y_{j,k}+z_k \geq 0, \quad \forall i, j, k \in [n].
% \end{align*}

{
\subsubsection{Union of Ellipsoidal Constraints}
We now demonstrate that an even broader class of problems than MIOs with logical constraints can be cast within our framework. Concretely, we demonstrate that constraints
$\bm{x} \in \mathcal{S}:=\bigcup_{i=1}^k (Q_i \cap P_i),$ where $Q_i:=\{\bm{x} \in \mathbb{R}^n: \bm{x}^\top \bm{Q}_i\bm{x}+\bm{h}_i^\top \bm{x}+g_i \leq 0\}$, with $\bm{Q}_i \succeq \bm{0}$, {is an ellipsoid} and  $P_i:=\{\bm{x}: \bm{A}_i\bm{x} \leq \bm{b}_i\}$ is a %bounded
polytope, can be reformulated as a special case of our framework. We remark that the constraint $\bm{x} \in \mathcal{S}$ is very general. Indeed, if we were to omit the quadratic constraints then we obtain a so-called {ideal} union of polyhedra formulation, which essentially all mixed-binary linear feasible regions admit \citep[see][]{vielma2015mixed}.

To derive a mixed-integer formulation with logical constraints of $\mathcal{S}$ that fits within our framework, we introduce $\bm{x}_i \in \mathbb{R}^n$ and  $\delta_i \in \{0, 1\}^n$, such that $\bm{x}_i \in Q_i \cap P_i$ if $\delta_i=1$, $\bm{x}_i = \bm{0}$ otherwise, and $\bm{x} = \sum_i \bm{x}_i$. We enforce $\bm{x}_i \in Q_i \cap P_i$ by introducing slack variables $\bm{\xi_i}$, $\rho_i$ for the linear and quadratic constraints respectively, and forcing them to be zero whenever $\delta_i =1$. Formally, $\mathcal{S}$ admits the following formulation
\begin{align}
    \bm{x}=\sum_{i=1}^k \bm{x}_i,\  \sum_{i=1}^k \delta_i=1,\\
    \bm{A}_i \bm{x}_i \leq \bm{b}_i + \bm{\xi}_i\quad \forall i \in [k], \nonumber\\
    \bm{x}_i^\top \bm{Q}_i \bm{x}_i+\bm{h}_i^\top \bm{x}_i + g_i \leq \rho_i\quad \forall i \in [k],\nonumber\\
    \bm{x}_i=\bm{0} \ \text{if} \ \delta_i=0\quad \forall i \in [k],\nonumber\\
    \bm{\xi}_i = 0 \ \text{if} \ (1-\delta_i)=0 \quad \forall i \in [k],\nonumber\\
    \rho_i = 0 \ \text{if} \ (1-\delta_i)=0 \quad \forall i \in [k].\nonumber
\end{align}

}

\subsection{A Regularization Assumption}
When we stated Problem \eqref{eqn:original_minlp}, we assumed that its objective function consists of a convex function $g(\bm{x})$ plus a regularization term $\Omega(\bm{x})$. We now formalize this assumption:
\begin{assumption} \label{assump:sep}
In Problem \eqref{eqn:original_minlp}, the regularization term $\Omega(\bm{x})$ is one of:
\begin{itemize}
\item a big-$M$ penalty function, $\Omega(\bm{x})= 0$ if $\ \Vert \bm{x} \Vert_\infty \leq M$ and $\infty$ otherwise, %which is equivalent to the big-$M$ formulation trick $ \| x_i \| \leq M z_i$ for each index $i$,
\item a ridge penalty, $ \Omega(\bm{x}) = \dfrac{1}{2\gamma} \| \bm{x} \|_2^2$.
\end{itemize}
\end{assumption}
This decomposition often constitutes a modeling choice in itself. We now illustrate this idea via the network design example.\FloatBarrier
\begin{box_example}\textit{
In the network design example \eqref{eqn:nd}, given the flow conservation structure $\bm{A} \bm{x} = \bm{b}$, we have that $\bm{x} \leq M \bm{e}$, where $M = \sum_{i: b_i > 0} b_i$. In addition, if $\bm{Q} \succ \mathbf{0}$ then the objective function naturally contains a ridge regularization term with $1/\gamma$ equal to the smallest eigenvalue of $\bm{Q}$.
Moreover, it is possible to obtain a tighter natural ridge regularization term by solving the following auxiliary semidefinite optimization problem a priori
\begin{align*}
\max_{\bm{q} \geq \bm{0}}\  \bm{e}^\top \bm{q} \quad & \text{s.t.} \quad \bm{Q}-\mathrm{Diag}(\bm{q}) \succeq \bm{0},
\end{align*}
and using $q_i$ as the ridge regularizer for each index $i$ \citep{frangioni2007sdp}.}
\end{box_example}
%\FloatBarrier
\vspace{5mm}
Big-$M$ constraints are often considered to be a modeling trick. However, our framework demonstrates that imposing either big-$M$ constraints or a ridge penalty is a regularization method, rather than a modeling trick. Interestingly, ridge regularization accounts for the relationship between the binary and continuous variables just as well as big-$M$ regularization, without performing an algebraic reformulation of the logical constraints\footnote{Specifically, ridge regularization enforces logical constraints through perspective functions, as is made clear in Section \ref{sec:perspectivecuts}.}.

Conceptually, both regularization functions are equivalent to a soft or hard constraint on the continuous variables $\bm{x}$. However, they admit practical differences:
{ For big-$M$ regularization, there usually
%\footnote{But not always; see \cite{hijazi2016constraint} for a (somewhat pathological) counterexample where no such $M_0$ exists.}
exists a finite value ${M}_0$, typically unknown a priori, such that if $M < {M}_0$, the regularized problem is infeasible. Alternatively, for every value of the ridge regularization parameter $\gamma$, if the original problem is feasible then the regularized problem is also feasible. Consequently, if there is no natural choice of $M$ then imposing ridge regularization may be less restrictive than imposing big-$M$ regularization.
However, for any $\gamma > 0$, the objective of the optimization problem with ridge regularization is different from its unregularized limit as $\gamma \rightarrow \infty$, while for big-$M$ regularization, there usually exists a finite value $M_1$ above which the two objective values match.}
%It could be argued that this addition is artificial and might result in some loss of generality. However, it is usually not possible to improve beyond complete enumeration of all feasible choices of $\bm{z}$ without making such an assumption. For instance, existing branch-and-bound approaches require a big-$M$ assumption to obtain lower bounds via continuous relaxations, at least whenever no ``natural'' value for $M$ exists. Consequently, we do not believe that our regularization assumption is too restrictive.
{We illustrate this discussion numerically in Section \ref{ssec:UCnumres}.}
\subsection{Duality to the Rescue}
In this section, we derive Problem \eqref{eqn:subproblem}'s dual and reformulate $f(\bm{z})$ as a maximization problem.
This reformulation is significant for two reasons: First, as shown in the proof of Theorem \ref{thm:innermax}, it leverages a non-linear reformulation of the logical constraints ``$x_i = 0$ if $z_i = 0$'' by introducing additional variables $v_i$ such that $v_i = z_i x_i$. Second, it proves that the regularization term $\Omega(\bm{x})$ drives the convexity and smoothness of $f(\bm{z})$, and thereby drives the computational tractability of the problem. To derive Problem
\eqref{eqn:subproblem}'s dual, we require:% the following assumption:
\begin{assumption}\label{cq}
For each subproblem generated by $f(\bm{z})$, where $\bm{z} \in \mathcal{Z}$, either the optimization problem is infeasible, or strong duality holds.
\end{assumption}

Note that all seven problems stated in Section \ref{sec:examples} satisfy Assumption \ref{cq}, as their inner problems are convex quadratics with linear or semidefinite constraints \citep[][Section 5.2.3]{boyd2004convex}. %Moreover, if we augment these problems with non-linear constraints then Assumption \ref{cq} holds under a constraint qualification, which can be verified by solving a separation problem \citep[see][]{ben2001lectures}.
Under Assumption \ref{cq}, the following theorem reformulates Problem \eqref{eqn:masterproblem} as a saddle-point problem:
%derives a saddle-point formulation of Problem \eqref{eqn:masterproblem}:
\begin{theorem}
\label{thm:innermax}
Under Assumption \ref{cq}, Problem \eqref{eqn:masterproblem} is equivalent to:% the following problem:
\begin{equation}
\begin{aligned}\label{eqn:saddlepointproblem}
    \min_{\bm{z} \in \mathcal{Z}} \ \max_{\bm{\alpha} \in \mathbb{R}^n} \quad & \bm{c}^\top \bm{z} + h(\bm{\alpha}) - \sum_{i=1}^n z_i \, \Omega^\star ( {\alpha}_i ),
\end{aligned}
\end{equation}
where $h(\bm{\alpha}) := \inf_{\bm{v}} \: g(\bm{v}) - \bm{v}^\top \bm{\alpha}$ is, up to a sign, the Fenchel conjugate of $g$ \citep[see][Chap. 3.3]{boyd2004convex}, and
\begin{equation*}
\begin{aligned}
\Omega^\star ( \beta ) &:= M \vert \beta \vert &\text{ for the big-$M$ penalty,} \\
\Omega^\star ( \beta ) &:= \tfrac{\gamma}{2} \beta^2 &\text{ for the ridge penalty.}
\end{aligned}
\end{equation*}
\end{theorem}
\proof
Let us fix some $\bm{z} \in \{0,1\}^n$, and suppose that strong duality holds for the inner minimization problem which defines $f(\bm{z})$. Then, after introducing additional variables $\bm{v} \in \mathbb{R}^n$ such that $v_i = z_i x_i$, we have
\begin{align*}
   f(\bm{z}) &= \min_{\bm{x},\bm{v}} \: g(\bm{v}) +\Omega(\bm{x}) \quad \mbox{  s.t.  } \bm{v} = \mathrm{Diag}(\bm{z})\bm{x}.
\end{align*}
Let $\bm{\alpha}$ denote the dual variables associated with the coupling constraint $\bm{v} = \mathrm{Diag}(\bm{z})\bm{x}$. The minimization problem is then equivalent to its dual problem, which is given by:
\begin{align*}
   f(\bm{z}) &= \max_{\bm{\alpha}} \: h(\bm{\alpha}) + \min_{\bm{x}} \left[ \Omega(\bm{x}) + \bm{\alpha}^\top \mathrm{Diag}(\bm{z}) \bm{x} \right],
\end{align*}
Since $\Omega(\cdot)$ is decomposable, i.e., $\Omega(\bm{x})=\sum_i \Omega_i({x}_i)$, we obtain:
\begin{align*}
\min_{\bm{x}} \left[ \Omega(\bm{x}) + \bm{\alpha}^\top \mathrm{Diag}(\bm{z}) \bm{x} \right]
&= \sum_{i=1}^n \min_{{x}_i} \left[ \Omega_i({x}_i) + {z}_i {x}_i {\alpha}_i \right]\\
&= \sum_{i=1}^n  - \Omega^\star(- {z}_i {\alpha}_i)=  - \sum_{i=1}^n {z}_i \Omega^\star({\alpha}_i),
\end{align*}
where the last equality holds as ${z}_i >0$ for the big-$M$ and ${z}_i^2 = {z}_i$ for the ridge penalty.

Alternatively, if the inner minimization problem defining $f(\bm{z})$ is infeasible, then its dual problem is unbounded by weak duality\footnote{Weak duality implies that the dual problem is either unfeasible or unbounded. Since the feasible set of the maximization problem does not depend on $\bm{z}$, it is always feasible, unless the original problem \eqref{eqn:original_minlp} is itself infeasible. Therefore, we assume without loss of generality that it is unbounded.}.
\endproof

\begin{remark}
%As pointed out by a referee of a previous version of this paper,
Without regularization, i.e., $\Omega(\bm{x})=0$, {a similar proof shows that} Problem \eqref{eqn:masterproblem} admits an interesting saddle-point formulation:
\begin{align*}
    \min_{\bm{z} \in \mathcal{Z}} \ \max_{\bm{\alpha} \in \mathbb{R}^n} \quad & \bm{c}^\top \bm{z} + h(\bm{\alpha}) \ \text{s.t.} \ \alpha_i=0, \ \text{if}\ z_i=1\quad \forall i \in [n],
\end{align*}
since $\Omega^\star(\alpha)= {\min_{{x}}  \left[ {x} {\alpha} -  \Omega({x}) \right] =} 0$ if $\alpha=0$, and $+\infty$ otherwise.
Consequently, the regularized formulation can be regarded as a relaxation of the original problem where the hard constraint $\alpha_i=0 \ \text{if} \ z_i=1$ is replaced with a soft penalty term $-z_i \Omega^\star(\alpha_i)$.
\end{remark}

\begin{remark}
The proof of Theorem \ref{thm:innermax} exploits three attributes of the regularizer $\Omega(\bm{x})$. Namely, (1) decomposability, i.e., $\Omega(\bm{x})=\sum_i \Omega_i(x_i)$, for appropriate scalar functions $\Omega_i$, (2) the convexity of $\Omega(\bm{x})$ in $\bm{x}$, and (3) the fact that $\Omega(\cdot)$ regularizes\footnote{Importantly, the third attribute allows us to strengthen the formulation by not associating $\bm{z}$ with $\bm{x}$ in $\Omega(\bm{x})$, since $\bm{x}_i=0$ is a feasible, indeed optimal choice of $\bm{x}$ for minimizing the regularizer when $z_i=0$; this issue is explored in more detail in \citep[Lemma 1]{bertsimas2018portfolio}; see also \citep[Appendix A.1]{bertsimas2019cosice}.} towards $0$, i.e., $\bm{0} \in \arg\min_{\bm{x}} \Omega(\bm{x})$. However, the proof does not explicitly require that $\Omega(\bm{x})$ is either a big-$M$ or a ridge regularizer. This suggests that our framework could be extended to other regularization functions. %, such as an $\ell_1$ (lasso) regularizer. We have refrained from incorporating a lasso regularizer in our numerical experiments, as its Fenchel conjugate, the $\ell_\infty$ norm, does not give rise to computationally useful subgradients.
\end{remark}

\begin{box_example}\textit{
For the network design problem \eqref{eqn:nd}, we have
\begin{align*}
h(\bm{\alpha}) &= \min_{\bm{x} \geq \mathbf{0} : \bm{A} \bm{x} = \bm{b}} \quad \tfrac{1}{2} \bm{x}^\top \bm{Q} \bm{x} + (\bm{d} - \bm{\alpha})^\top \bm{x}, \\
&= \max_{\bm{\beta}_0 \geq \bm{0}, \bm{p}} \quad \bm{b}^\top \bm{p} - \tfrac{1}{2} \left( \bm{A}^\top \bm{p} - \bm{d} + \bm{\alpha} + \bm{\beta}_0 \right)^\top \bm{Q}^{-1} \left( \bm{A}^\top \bm{p} - \bm{d} + \bm{\alpha} + \bm{\beta}_0 \right).
\end{align*}
Introducing $\boldsymbol{\xi} = \bm{Q}^{-1/2} \left( \bm{A}^\top \bm{p} - \bm{d} + \boldsymbol{\alpha} + \bm{\beta}_0 \right)$, we can further write
\begin{align*}
h(\bm{\alpha}) &= \max_{\boldsymbol{\xi}, \bm{p}} \quad \bm{b}^\top \bm{p} - \tfrac{1}{2} \| \boldsymbol{\xi}\|_2^2 \: \mbox{ s.t } \: \bm{Q}^{1/2}\boldsymbol{\xi} \geq \bm{A}^\top \bm{p} - \bm{d} + \boldsymbol{\alpha}.
\end{align*}
Hence, Problem \eqref{eqn:nd} is equivalent to minimizing over $\bm{z} \in \mathcal{Z}$ the function
\begin{align*}
\bm{c}^\top \bm{z} + f(\bm{z}) = \max_{\alpha, \boldsymbol{\xi}, \bm{p}} \quad & \bm{c}^\top \bm{z} + \bm{b}^\top \bm{p} - \tfrac{1}{2} \| \boldsymbol{\xi} \|_2^2  - \sum_{j=1}^n z_j \, \Omega^\star ( {\alpha}_j ) \\
\mbox{ s.t } \quad & \bm{Q}^{1/2}\boldsymbol{\xi} \geq \bm{A}^\top \bm{p} - \bm{d} + \boldsymbol{\alpha}.
\end{align*}}
\end{box_example}

Theorem \ref{thm:innermax} reformulates $f(\bm{z})$ as an inner maximization problem, namely
\begin{align} \label{eqn:innermax}
    f(\bm{z}) = \max_{\bm{\alpha} \in \mathbb{R}^n} \quad & h(\bm{\alpha}) - \sum_{i=1}^n z_i \, \Omega^\star ( {\alpha}_i ),
\end{align}
for any feasible binary $\bm{z} \in \mathcal{Z}$. The regularization term $\Omega$ will be instrumental in our numerical strategy for it directly controls both the convexity and smoothness of $f$. {Note that \eqref{eqn:innermax} extends the definition of $f(\bm{z})$ to the convex set $\mathrm{Bool}(\mathcal{Z})$, obtained by relaxing the constraints $\bm{z} \in \{0,1\}^p$ to $\bm{z} \in [0,1]^p$ in the definition of $\mathcal{Z}$.}
\subsubsection*{Convexity} $f(\bm{z})$ is convex in $\bm{z}$ as a point-wise maximum of linear function of $\bm{z}$. In addition, denoting $\bm{\alpha}^\star(\bm{z})$ a solution of \eqref{eqn:innermax}, we have the lower-approximation:
\begin{equation}
\label{eqn:linearlowerapp2}
    f(\Tilde{\bm{z}}) \geq f({\bm{z}}) + \nabla f({\bm{z}})^\top (\Tilde{\bm{z}} - \bm{z}) \quad \forall \Tilde{\bm{z}} \in \mathcal{Z},
\end{equation}
where $[\nabla f({\bm{z}})]_i := - \Omega^\star(\alpha^\star(\bm{z})_i)$ is a sub-gradient of $f$ at $\bm{z}$.

We remark that if the maximization problem in $\bm{\alpha}$ defined by $f(\bm{z})$ admits multiple optimal solutions then the corresponding lower-approximation of $f$ at $\bm{z}$ may not be unique. This behavior can severely hinder the convergence of outer-approximation schemes such as {Benders'} decomposition. Since the work of \cite{magnanti1981accelerating} on Pareto optimal cuts, many strategies have been proposed to improve the cut selection process in the presence of degeneracy \citep[see][Section 4.4 for a review]{fischetti2016benders}. However, the use of ridge regularization ensures that the objective function in \eqref{eqn:saddlepointproblem} is strongly concave in ${\alpha}_i$ such that $z_i > 0$, and therefore guarantees that there is a unique optimal choice of ${\alpha}_i^\star(\bm{z})$. In other words, ridge regularization naturally inhibits degeneracy.

\subsubsection*{Smoothness} $f(\bm{z})$ is smooth, in the sense of Lipschitz continuity, which is a crucial property for deriving bounds on the integrality gap of the Boolean relaxation, and designing local search heuristics in Section \ref{sec:algo}. Formally, the following proposition follows from Theorem \ref{thm:innermax}:% (proof omitted):
\begin{proposition} \label{thm:lipschitz}
For any $\bm{z},\bm{z}' \in \mathrm{Bool}\left(\mathcal{Z}\right)$,
\begin{enumerate}
\item[(a)] With big-$M$ regularization,
% \begin{align*}
      $f(\bm{z}')-f(\bm{z})  \leq M  \displaystyle \sum_{i=1}^n (z_i - z'_i) | \alpha^\star(\bm{z}')_i |$.
% \end{align*}
\item[(b)] With ridge regularization,
% \begin{align*}
      $f(\bm{z}')-f(\bm{z})  \leq \dfrac{\gamma}{2} \displaystyle \sum_{i=1}^n (z_i - z'_i) \alpha^\star(\bm{z}')_i^2$.
% \end{align*}
\end{enumerate}
\end{proposition}
\proof
By Equation \eqref{eqn:saddlepointproblem},
\begin{align*}
f(\bm{z}')-f(\bm{z}) & = \max_{\bm{\alpha}' \in \mathbb{R}^n} \left(h(\bm{\alpha}')-\sum_{i=1}^n z'_i \Omega^\star(\alpha'_i)\right)-\max_{\bm{\alpha} \in \mathbb{R}^n} \left(h(\bm{\alpha})-\sum_{i=1}^n z_i \Omega^\star(\alpha_i)\right),\\
& = h(\bm{\alpha}^\star(\bm{z}'))-\sum_{i=1}^n z'_i \Omega^\star(\alpha^\star(\bm{z}')_i) - h(\bm{\alpha}^\star(\bm{z}'))+\sum_{i=1}^n z_i \Omega^\star(\alpha^\star(\bm{z}')_i),\\
& \leq \sum_{i=1}^n (z_i-z_i')     \Omega^\star(\alpha^\star(\bm{z}')_i),
\end{align*}
where the inequality holds because {an} optimal choice of $\bm{\alpha}'$ is a feasible choice of $\bm{\alpha}$.
\endproof

Proposition \ref{thm:lipschitz} demonstrates that, when the coordinates of $\bm{\alpha}^\star (\bm{z})$ are uniformly bounded\footnote{Such a uniform bound always exists, as $f(\bm{z})$ is only supported on a finite number of binary points. Moreover, the strong concavity of $h$ can yield stronger bounds (see Appendix \ref{sec:lipcont}).} with respect to $\bm{z}$, $f(\bm{z})$ is Lipschitz-continuous, with a constant {$L$} proportional to $M$ (resp. $\gamma$) in the big-$M$ (resp. ridge) case. { We provide explicit bounds on the magnitude of $L$ in Appendix \ref{sec:lipcont}.}

\subsection{Merits of Ridge, Big-\texorpdfstring{$M$}{M} Regularization: Theoretical Perspective}\label{ssec:meritstheory}
{
In this section, we {propose} a framework to reformulate MINLOs with logical constraints, which comprises regularizing MINLOs via either the widely used big-$M$ modeling paradigm or the less popular ridge regularization paradigm. We summarize the advantages and disadvantages of each regularizer in Table \ref{tab:comparison_regularization}. However, note that we have not yet established how these characteristics impact the numerical tractability and quality of the returned solution; this is the topic of the next two sections.
\begin{table}[h]
\centering
\caption{Summary of the advantages ($+$) /disadvantages ($-$) of both techniques.}
\begin{tabular}{ll} \toprule
Regularization & Characteristics \\ \midrule
\multirow{3}{*}{Big-$M$} & ($+$) Linear constraints\\
& ($+$) Supplies the same objective if $M > M_1$, for some $M_1 < \infty$ \\
& ($-$) Leads to infeasible problem if $M < M_0$, for some $M_0 < \infty$ \\
\midrule
\multirow{3}{*}{Ridge} & ($+$) Strongly convex objective\\
& ($-$) Systematically leads to a different objective for any $\gamma > 0$ \\
& ($+$) Preserves the feasible set\\
\bottomrule
\end{tabular}
\label{tab:comparison_regularization}
\end{table}
}

\section{An Efficient Numerical Approach} \label{sec:algo}

We now present an efficient numerical approach to solve Problem \eqref{eqn:saddlepointproblem}. The backbone is an outer-approximation strategy, embedded within a branch-and-bound procedure to solve the problem exactly. We also propose local search and rounding heuristics to find good feasible solutions, and use information from the Boolean relaxation to improve the duality gap.

\subsection{Overall Outer-Approximation Scheme} \label{sec:algo.oa}
Theorem \ref{thm:innermax} reformulates the function $f(\bm{z})$ as an inner maximization problem, and demonstrates that $f(\bm{z})$ is convex in $\bm{z}$, meaning a linear outer approximation provides a valid underestimator of $f(\bm{z})$, as outlined in Equation \eqref{eqn:linearlowerapp2}.
% \begin{equation}
% \label{eqn:linearlowerapp}
%     f(\Tilde{\bm{z}}) \geq f({\bm{z}}) + \nabla f({\bm{z}})^\top (\Tilde{\bm{z}} - \bm{z}), \, \forall \Tilde{\bm{z}},
% \end{equation}
% where $[\nabla f({\bm{z}})]_i := - \Omega^\star(\alpha^\star(\bm{z})_i)$ is a sub-gradient of $f$ at $\bm{z}$.
Consequently, a valid numerical strategy for minimizing $f(\bm{z})$ is to iteratively minimize a piecewise linear lower-approximation of $f$ and refining this approximation at each step until some approximation error $\varepsilon$ is reached, as described in Algorithm \ref{alg:OA}. This scheme was originally proposed for continuous decision variables by \cite{kelley1960cutting},
% independently proposed for continuous decision variables by \cite{cheney1959newton} and \cite{kelley1960cutting},
and later extended to binary decision variables by \cite{duran1986outer}, who provide a proof of termination in a finite, yet exponential in the worst case, number of iterations.
\begin{algorithm*}
\caption{Outer-approximation scheme}
\label{alg:OA}
\begin{algorithmic}
\REQUIRE Initial solution $\bm{z}^1$
\STATE $t \leftarrow 1 $
\REPEAT
\STATE Compute $\bm{z}^{t+1}, \eta^{t+1}$ solution of
\begin{align*}
\min_{\bm{z} \in \mathcal{Z}, \eta} \: \bm{c}^\top \bm{z} + \eta \quad \mbox{ s.t. } \forall s \in \{ 1, \dots, t\}, \: \eta \geq f(\bm{z}^s) + \nabla f (\bm{z}^s)^\top (\bm{z}-\bm{z}^s)
\end{align*}
\STATE Compute $f(\bm{z}^{t+1})$ and $\nabla f (\bm{z}^{t+1})$
\STATE $t \leftarrow t+1 $
\UNTIL{$f(\bm{z}^{t+1}) - \eta^{t+1} \leq \varepsilon$}
\RETURN $\bm{z}^t$
\end{algorithmic}
\end{algorithm*}

To avoid solving a mixed-integer linear optimization problem at each iteration, as suggested in the pseudo-code, this strategy can be integrated within a single branch-and-bound procedure using lazy callbacks, as originally proposed by \cite{quesada1992lp}. Lazy callbacks are now standard tools in commercial solvers such as Gurobi and CPLEX and provide significant speed-ups for outer-approximation algorithms. With this implementation, the commercial solver constructs a single branch-and-bound tree and generates  {a new cut at} a feasible solution $\bm{z}$.

We remark that the second-stage minimization problem may be infeasible at some $\bm{z}^t$. In this case, we generate a feasibility cut rather than outer-approximation cut. In particular, the constraint $\sum_{i}  z^t_i (1-z_i) + \sum_{i}(1-z^t_i) z_i \geq 1$ excludes the iterate $\bm{z}^t$ from the feasible set. {Stronger feasibility cuts can be obtained by leveraging problem specific structure. For instance, when the feasible set satisfies $\bm{z}^t \notin \mathcal{Z} \implies \forall \bm{z} \leq \bm{z}^t,\ \bm{z} \notin \mathcal{Z}$, $\sum_{i}(1-z^t_i) z_i \geq 1$ is a valid feasibility cut. % (b) if the second-stage problem is a linear optimization problem, an extreme ray with positive marginal cost defines a feasibility cut.}.
Alternatively, one can %An alternative approach for generating feasibility cuts is to
invoke conic duality if %This approach involves assuming that
$g(\bm{x})$ generates a conic feasibility problem. Formally, assume $$g(\bm{x})=\begin{cases} \langle \bm{c}, \bm{x} \rangle, & \text{if} \ \bm{A}\bm{x}=\bm{b},\ \bm{x} \in \mathcal{K},\\ +\infty, & \text{otherwise,} \end{cases}$$ where $\mathcal{K}$ is a closed convex cone.
This assumption gives rise to some loss of generality. Note, however, that all the examples in the previous section admit conic reformulations by taking appropriate Cartesian products of the linear, second-order and semidefinite cones \citep[][]{ben2001lectures}. Assuming that $g(\bm{x})$ is of the prescribed form, we have the dual conjugate
\begin{align*}
    h(\bm{\alpha})=\inf_{\bm{x}}\langle \bm{x}, \bm{\alpha}\rangle-g(\bm{x})=\max_{\bm{\pi}} \langle \bm{b}, \bm{\pi}\rangle +\begin{cases}
    0, & \text{if} \ \bm{c}-\bm{\alpha}-\bm{A}^\top \bm{\pi} \in \mathcal{K}^\star,\\
    +\infty, & \text{otherwise},
    \end{cases}
\end{align*}
where $\mathcal{K}^\star$ is the dual cone to $\mathcal{K}$. In this case, if some binary vector $\bm{z}$ gives rise to an infeasible subproblem, i.e., $f(\bm{z})=+\infty$, then the conic duality theorem implies\footnote{We should note that this statement is, strictly speaking, not true unless we impose regularization. Indeed, the full conic duality theorem \citep[Theorem 2.4.1]{ben2001lectures} allows for the possibility that a problem is infeasible but asymptotically feasible, i.e., $$\nexists \bm{x}: \bm{A}\bm{x}=\bm{b}, \bm{x} \in \mathcal{K} \ \text{but} \ \exists \{ \bm{x}_t\}_{t=1}^\infty: \bm{x}_t \in \mathcal{K}\ \forall t \ \text{with}\  \Vert \bm{A}\bm{x}_t -\bm{b}\Vert \rightarrow 0.$$ Fortunately, the regularizer $\Omega(\bm{x})$ alleviates this issue, because it is coercive (i.e., ``blows up'' to $+\infty$ as $\Vert\bm{x}\Vert \rightarrow \infty$) and therefore renders all unbounded solutions infeasible and ensures the compactness of the level sets of $g(\bm{x})+\Omega(\bm{x})$.} that there is a \textit{certificate} of infeasibility $(\bm{\alpha}, \bm{\pi}$) such that
\begin{align*}
    \bm{c}-\bm{\alpha}-\bm{A}^\top \bm{\pi} \in \mathcal{K}^\star, \langle \bm{b}, \bm{\pi} \rangle >\sum_{i=1}^n z_i \Omega^\star(\alpha_i).
\end{align*}
Therefore, to restore feasibility, we can simply impose the cut
$
    \langle \bm{b}, \bm{\pi} \rangle \leq \sum_{i=1}^n z_i \Omega^\star(\alpha_i).
$
}

As mentioned in Section \ref{ssec:background}, the rate of convergence of outer-approximation schemes depends heavily on three criterion. We now provide practical guidelines on how to meet these criterion:
\begin{enumerate}
    \item \emph{Fast cut generation strategy: } To generate a cut, one solves the second-stage minimization problem \eqref{eqn:subproblem} (or its dual) in $\bm{x}$, which contains no discrete variables and is usually orders of magnitude faster to solve than the original mixed-integer problem \eqref{eqn:original_minlp}. Moreover, the minimization problem in $\bm{x}$ needs to be solved only for the coordinates $x_i$ such that $z_i=1$. In practice, this approach yields a sequence of subproblems of much smaller size than the original problem, especially if $\mathcal{Z}$ contains a cardinality constraint. For instance, for the sparse empirical risk minimization problem \eqref{eqn:serm}, each cut is generated by solving a subproblem with $n$ observations and $k$ features, where $k \ll p$. For this reason, we recommend generating cuts at binary $\bm{z}$'s, which are often sparser than continuous $\bm{z}$'s. This recommendation can be relaxed in cases where the separation problem can be solved efficiently even for dense $\bm{z}$'s; for instance, in uncapacitated facility location problems, each subproblem is a knapsack problem which can be solved by sorting \citep{fischetti2016redesigning}. {If possible, we recommend theoretically analyzing the sparsity of the optimal solution a priori, to derive an explicit} cardinality or budget constraint on $\bm{z}$ and ensure the sparsity of each incumbent solution.
    \item \emph{Cut selection rule in presence of degeneracy: } In the presence of degeneracy, selection criteria, such as Pareto optimality \citep{magnanti1981accelerating}, have been proposed to accelerate convergence. However, these criteria are numerous, computationally expensive and all in all, can do more harm than good \citep{papadakos2008practical}.
    %These solutions have merit, but considered the problem formulation as given and immutable.
    In an opposite direction, we recommend alleviating the burden of degeneracy by design, by imposing a ridge regularizer whenever degeneracy hinders convergence.
    \item \emph{Rich root node analysis: } As suggested in \cite{fischetti2016benders}, providing the solver with as much information as possible at the root node can drastically improve convergence of outer-approximation methods. This is the topic of the next two sections. Restarting mechanisms, as described in \cite[][Section 5.2]{fischetti2016benders}, could also be useful, although we do not implement them in the present paper.%any of those in the present paper.
\end{enumerate}
{These ingredients, and especially the ability to generate cuts efficiently, dictate which types of problems could benefit the most from our approach and which regularizer to use. Problems with an explicit cardinality constraint, for instance, would require a small subproblem to be solved at each iteration. For network design problems, the network flow structure of the feasible set is a key numerical asset so we intuit that ridge regularization, which leaves the feasible set unchanged, would be very efficient. On the other hand, for uncapacitated facility location, sub-problems with big-$M$ regularization boils down to a knapsack problem and can be solved efficiently via sorting, as discussed in \citep[Section 3.1]{fischetti2016redesigning}. }

\subsection{Improving the Lower-Bound: A Boolean Relaxation} \label{ssec:boolrelax}
To certify optimality, high-quality lower bounds are of interest and can be obtained by relaxing the integrality constraint $\bm{z} \in \{0,1\}^n$ in the definition of $\mathcal{Z}$ to $\bm{z} \in [0,1]^n$. In this case, the Boolean relaxation of \eqref{eqn:masterproblem} is:
\begin{align*}
    \min_{\bm{z} \in \mathrm{Bool}({\mathcal{Z}})} \quad & \bm{c}^\top \bm{z} + f(\bm{z}),
\end{align*}
which can be solved using Kelley's algorithm \citep{kelley1960cutting}, which is a continuous analog of Algorithm \ref{alg:OA}. Stabilization strategies have been empirically successful to accelerate the convergence of Kelley's algorithm, as recently demonstrated on uncapacitated facility location problems by \cite{fischetti2016redesigning}. However, for Boolean relaxations, Kelley's algorithm computes $f(\bm{z})$ and $\nabla f (\bm{z})$ at dense vectors $\bm{z}$, which is (sometimes substantially) more expensive than for sparse binary vectors $\bm{z}$'s, unless each subproblem can be solved efficiently as in \cite{fischetti2016redesigning}.

Alternatively, the continuous minimization problem admits a reformulation
\begin{equation}
\label{eqn:relaxation_minmax}
    \min_{\bm{z} \in \mathrm{Bool}(\mathcal{Z}) } \: \max_{\bm{\alpha} \in \mathbb{R}^m} \quad \bm{c}^\top \bm{z} + h(\bm{\alpha}) - \sum_{i=1}^n z_i \, \Omega^\star ( {\alpha}_j ).
\end{equation}
analogous to Problem \eqref{eqn:saddlepointproblem}. Under Assumption \ref{cq}, we can further write the min-max relaxation formulation \eqref{eqn:relaxation_minmax} as a non-smooth maximization problem
\begin{align*}
    \max_{\bm{\alpha} \in \mathbb{R}^n} \: q(\bm{\alpha}), \quad \mbox{ with } \quad q(\bm{\alpha}) := h(\bm{\alpha}) + \min_{\bm{z} \in \mathrm{Bool}(\mathcal{Z})} \: \sum_{i=1}^n \left( c_i - \Omega^\star(\alpha_i) \right) z_i
\end{align*}
and apply a projected sub-gradient ascent method as in \cite{bertsimas2019sparse}. We refer to \cite[][Chapter 7.5.]{bertsekas1997nonlinear} for a discussion on implementation choices regarding step-size schedule and stopping criteria, and \cite{renegar2018simple} for recent enhancements using restarting.

The benefit from solving the Boolean relaxation with these algorithms is threefold. First, it provides a lower bound on the objective value of the discrete optimization problem \eqref{eqn:masterproblem}. Second, it generates valid linear lower approximations of $f(\bm{z})$ to initiate the cutting-plane algorithm with. Finally, it supplies a sequence of continuous solutions that can be rounded and polished to obtain good binary solutions. Indeed, the Lipschitz continuity of $f(\bm{z})$ suggests that high-quality feasible binary solutions can be found in the neighborhood of a solution to the Boolean relaxation. We formalize this observation in the following theorem:
\begin{theorem} \label{thm:integrality.gap}
Let $\bm{z}^\star$ denote a solution to the Boolean relaxation \eqref{eqn:relaxation_minmax}, $\mathcal{R}$ denote the indices of $\bm{z}^\star$ with fractional entries, and $\alpha^\star(\bm{z})$ denote a best choice of $\bm{\alpha}$ for a given $\bm{z}$. Suppose that for any $\bm{z} \in \mathcal{Z}$, $|\alpha^\star(\bm{z})_j| \leq L$. Then, a random rounding $\bm{z}$ of $\bm{z}^\star$, i.e., $z_j \sim Bernoulli(z^\star_j)$, satisfies
% \begin{align*}
    $0 \leq f(\bm{z})-f(\bm{z}^\star) \leq \epsilon$
% \end{align*}
with probability at least $p = 1- |\mathcal{R}| \exp \left( -\tfrac{\epsilon^2}{\kappa} \right),$ where
\begin{equation*}
\begin{aligned}
\kappa &:= 2 M^2 L^2 |\mathcal{R}|^2 &\text{ for the big-$M$ penalty,} \\
\kappa &:= \tfrac{1}{2} \gamma^2 L^4 |\mathcal{R}|^2  &\text{ for the ridge penalty.}
\end{aligned}
\end{equation*}
\end{theorem}
{ We provide a formal proof of this result in Appendix \ref{sec:proof.randomrounding}.}  This result calls for multiple remarks:
\begin{itemize}
\item For $\varepsilon > \sqrt{\kappa \ln(|\mathcal{R}|)}$, we have that $p > 0$, which implies the existence of a binary $\varepsilon$-optimal solution in the neighborhood of $\bm{z}^\star$, which in turn bounds the integrality gap by $\varepsilon$. As a result, lower values of $M$ or $\gamma$ typically make the discrete optimization problem easier.
\item A solution to the Boolean relaxation often includes some binary coordinates, i.e., $|\mathcal{R}| < n$. In this situation, it is tempting to fix $z_i = z^\star_i$ for $i \notin \mathcal{R}$ and solve the master problem \eqref{eqn:masterproblem} over coordinates in $\mathcal{R}$. In general, this approach provides sub-optimal solutions. However, Theorem \ref{thm:integrality.gap} quantifies the price of fixing variables and bounds the optimality gap by $\sqrt{\kappa \ln(|\mathcal{R}|)}$.
\item In the above high-probability bound, we do not account for the feasibility of the randomly rounded solution $\bm{z}$. Accounting for $\bm{z}$'s feasibility marginally reduces the probability given above, as shown for general discrete optimization problems by \cite{raghavan1987randomized}.
%\item If $\mathcal{R}$ is empty, by the probabilistic method, the relaxation is tight and $\bm{z}^\star$ solves Problem \eqref{eqn:saddlepointproblem}.
\end{itemize}

%We remark that the saddle-point formulation \eqref{eqn:relaxation_minmax} can be used to derive sufficient conditions for the Boolean relaxation to be tight; we omit these details in this condensed version of the manuscript.

Under specific problem structure, other strategies might be more efficient than Kelley's method or the subgradient algorithm. For instance, if $\mathrm{Bool}\left(\mathcal{Z}\right)$ is a polyhedron, then the inner minimization problem defining $q(\bm{\alpha})$ is a linear optimization problem that can be rewritten as a maximization problem by invoking strong duality. Although we only consider linear relaxations here, tighter bounds could be attained by taking a higher-level relaxation from a relaxation hierarchy, such as the \cite{lasserre2001explicit} hierarchy \citep[see][for a comparison]{laurent2003comparison}. {The main benefit of such a relaxation is that while the aforementioned Boolean relaxation only controls the first moment of the probability measure studied in Theorem \ref{thm:integrality.gap}, higher level relaxations control an increasing sequence of moments of the probability measure and thereby provide non-worsening probabilistic guarantees for randomized rounding methods.} However, the additional tightness of these bounds comes at the expense of solving relaxations with additional variables and constraints\footnote{$n^2$ additional variables and $n^2$ additional constraints for empirical risk minimization, versus $n+1$ additional variables and $n$ additional constraints for the linear relaxation.}; yielding a sequence of ever-larger semidefinite optimization problems. Indeed, even the SDP relaxation which controls the first two moments of a randomized rounding method is usually intractable when $n>300$, with current technology. For an analysis of higher-level relaxations in sparse regression problems, we refer the reader to \cite{atamturk2019rank}.

\subsection{Improving the Upper-Bound: Local Search and Rounding}
To improve the quality of the upper-bound, i.e., the cost associated with the best feasible solution found so far, we implement two rounding and local-search strategies.

Our first strategy is a randomized rounding strategy, which is inspired by Theorem \ref{thm:integrality.gap}. Given $\bm{z}_0 \in \mathrm{Bool}(\mathcal{Z})$, we generate randomly rounded vectors $\bm{z}$ {by sampling $\bm{z}$ according to $z_i \sim \mathrm{Bernoulli}(z_{0i})$ until $\bm{z} \in \mathcal{Z}$, which happens with high probability since $\mathbb{E}[\bm{z}] = \bm{z}_0$ satisfies all the constraints which describe $\mathcal{Z}$, besides integrality \citep{raghavan1987randomized}.}

Our second strategy is a sequential rounding procedure, which is informed by the lower-approximation on $f(\bm{z})$, as laid out in Equation \eqref{eqn:linearlowerapp2}. Observing that the $i$th coordinate $\nabla f(\bm{z}_0)_i$ provides a first-order indication of how a change in $z_i$ might impact the overall cost, we proceed in two steps. We first round down all coordinates such that $\nabla f(\bm{z}_0)_i(0 - \bm{z}_{0i}) < 0$. Once the linear approximation of $f$ only suggests rounding up, we round all coordinates of $\bm{z}$ to $1$ and iteratively bring some coordinates to $0$ to restore feasibility.

If $\bm{z}_0$ is binary, we implement a comparable local search strategy. If $\bm{z}_{0i} = 0$, then switching the $i$th coordinate to one increases the cost by at least $\nabla f(\bm{z}_0)_i$. Alternatively, if $\bm{z}_{0i} = 1$, then switching it to zero increases the cost by at least $-\nabla f(\bm{z}_0)_i$. We therefore compute the one-coordinate change which provides the largest potential cost improvement. However, as we only have access to a lower approximation of $f$, we are not guaranteed to generate a cost-decreasing sequence. Therefore, we terminate the procedure as soon as it cycles. A second complication is that, due to the constraints defining $\mathcal{Z}$,
the best change sometimes yields an infeasible $\bm{z}$. In practice, for simple constraints such as $\bm{\ell} \leq \bm{z} \leq \bm{u}$, we forbid switches which break feasibility; for cardinality constraints, we perform the best switch and then restore feasibility at minimal cost when necessary.

\subsection{Relationship With Perspective Cuts} \label{sec:perspectivecuts}
In this section, we connect the perspective cuts introduced by \cite{frangioni2006perspective} with our framework and discuss the merits of both approaches, in theory and in practice. To the best of our knowledge, a connection between Boolean relaxations of  the two approaches has only been made in the context of sparse regression, by \cite{xie2018ccp}. That is, the general connection we make here between the discrete problems, as well as their respective cut generating procedures, is novel.

We first demonstrate that imposing the ridge regularization term $\Omega(\bm{x})=\tfrac{1}{2 \gamma}\Vert \bm{x}\Vert_2^2$ naturally leads to the perspective formulation of \cite{frangioni2006perspective}:
\begin{theorem}
\label{thm:perspective}
Suppose that $\Omega(\bm{x})=\tfrac{1}{2 \gamma}\Vert \bm{x}\Vert_2^2$ and that Assumption \ref{cq} holds. Then, Problem \eqref{eqn:saddlepointproblem} is equivalent to the following optimization problem:
\begin{equation} \label{eqn:perspectiveformulation}
\min_{ \bm{z} \in \mathcal{Z}} \:
\min_{\bm{x} \in \mathbb{R}^{n}} \quad
\bm{c}^\top \bm{z}+ g(\bm{x})
+ \frac{1}{2\gamma} \sum_{i=1}^n \begin{dcases} \frac{x_i^2}{z_i}, & \text{if } z_i > 0, \\
0, & \text{if } z_i=0 \text{ and } x_i = 0, \\
\infty, & \text{otherwise}. \end{dcases}
\end{equation}
\end{theorem}

Theorem \ref{thm:perspective} follows from taking the dual of the inner-maximization problem in Problem \eqref{eqn:innermax}{; see Appendix \ref{ssec:proofofperspectivethm} for a formal proof}. Note that the equivalence stated in Theorem \ref{thm:perspective} also holds for $\bm{z} \in \mathrm{Bool}(\mathcal{Z})$. As previously observed in \cite[][]{ben2001lectures,akturk2009strong}, Problem \eqref{eqn:perspectiveformulation} can be formulated as a second-order cone problem (SOCP)
\begin{equation} \label{eqn:perspective_socp}
\begin{aligned}
    \min_{\bm{x} \in \mathbb{R}^{n}, \bm{z} \in \mathcal{Z}, \bm{\theta} \in \mathbb{R}^n} \quad  \bm{c}^\top \bm{z} + g(\bm{x})+\sum_{i=1}^n \theta_i \quad
    \text{s.t.}  \ \left\Vert \begin{pmatrix} \sqrt{\tfrac{2}{\gamma}}x_i \\ \theta_i-z_i \end{pmatrix}\right\Vert_2 \leq \theta_i+z_i \quad \forall i \in [n].
\end{aligned}
\end{equation}
and solved by linearizing the SOCP constraints into so-called perspective cuts, i.e.,
% , i.e.,
% \begin{align*}
%     \min_{\bm{x} \in \mathbb{R}^n, \bm{z} \in \mathcal{Z}, \bm{\theta} \in \mathbb{R}^n} \quad \bm{c}^\top \bm{x} +  g(\bm{x}) + \bm{e}^\top \bm{\theta} \quad
%     \text{s.t.} \quad  \theta_i \geq \tfrac{1}{2\gamma} \bar{x}_i(2x_i-\bar{x}_i z_i), \ \forall i \in [n], \forall \bar{x} \in  \bar{\mathcal{X}}.
% \end{align*}
% where the constraints
${\theta_i \geq \tfrac{1}{2\gamma} \bar{x}_i(2x_i-\bar{x}_i z_i), \forall \bar{x} \in  \bar{\mathcal{X}}}$, which have been extensively studied in the literature in the past fifteen years \citep{frangioni2006perspective,gunluk2010perspective,dong2015regularization, frangioni2016approximated, atamturk2019rank}. Observe that by separating Problem \eqref{eqn:perspectiveformulation} into master and subproblems, an outer approximation algorithm yields the same cut \eqref{eqn:linearlowerapp2} as in our scheme. {In this regard, our approach supplies a new and insightful derivation of the perspective cut approach. It is worth noting that our proposal can easily be implemented within a standard integer optimization solver such as CPLEX or Gurobi using callbacks, while existing implementations of the perspective cut approach have required tailored branch-and-bound procedures \citep[see, e.g.,][Section $3.1$]{frangioni2006perspective}.

}%However, while theoretically similar, there are subtle differences which make the computational performance of our proposal more attractive:
% \begin{itemize}
%     \item Our outer-approximation approach only includes cuts corresponding to optimal choices of $\bm{x}$ for a given incumbent solution $\bm{z}$. Alternatively, the perspective cut approach is a {outer-approximation scheme} \citep[see][]{duran1986outer} which generates the optimal cut with respect to $(\bm{z},\bm{x})$. Consequently, the perspective cut approach possesses weaker convergence properties.
%     \item Our outer-approximation scheme can easily be implemented within a standard integer optimization solver such as CPLEX or Gurobi using callbacks. Unfortunately, an {outer-approximation approach} for the perspective formulation requires a tailored branch-and-bound procedure \citep[see][Section $3.1$ for details]{frangioni2006perspective}. In this regard, our approach is more practical.
% \end{itemize}
{
\subsection{Merits of Ridge, Big-\texorpdfstring{$M$}{M} Regularization: Algorithmic Perspective}\label{ssec:algmerits}
We now summarize the relative merits of applying either ridge or big-$M$ regularization from an algorithmic perspective:
\begin{itemize}
    \item As noted in our randomized rounding guarantees in Section \ref{ssec:boolrelax}, the two regularization methods provide comparable bound gaps when $2M \approx \gamma L$, while if $2M \ll \gamma L$, big-$M$ regularization provides smaller gaps, and if $2M \gg \gamma L$, ridge regularization provides smaller gaps.
    \item For linear problems, ridge regularization limits dual degeneracy, while big-$M$ regularization does not. This benefit, however, has to be put in balance with the extra runtime and memory requirements needed for solving a quadratic, instead of linear, separation problem.
\end{itemize}
In summary, the benefits of applying either big-$M$ or ridge regularization are largely even and depend on the specific instance to be solved. In the next section, we perform a sequence of numerical experiments on the problems studied in Section \ref{sec:examples}, to provide empirical guidance on which regularization approach works best when.
}

\section{Numerical Experiments}\label{sec:numresults}
In this section, we evaluate our single-tree cutting-plane algorithm, implemented in Julia 1.0 %\citep{bezanson2017julia}
using CPLEX $12.8.0$ and the Julia package \verb|JuMP.jl| version $0.18.4$ \citep{dunning2017jump}. We compare our method against solving the natural big-$M$ or MISOCP formulations directly, using CPLEX $12.8.0$. All experiments were performed on
%the \textit{engaging} cluster, a high performance cluster at the \textit{Massachusetts Green High Performance Computing Centre} (MGHPCC) which comprises a set of
one Intel Xeon E$5-2690$ v4 $2.6$GHz CPU core and using $32$ GB RAM. %The big-$M$ and MISOCP formulations are allocated $8$ threads, unless explicitly stated otherwise. However, \verb|JuMP.jl| version $0.18.4$ does not allow our outer-approximation scheme to benefit from multi-threading.%, we allocate our cutting-plane methods a single thread. %is currently not thread-safe and CPLEX currently cannot combine multiple threads with lazy callbacks and non-thread-safe code, we allocate our cutting-plane methods a single thread.

% We perform those experiments on real-world and synthetic instances of all problems presented in Section \ref{sec:examples}, with the exception of the empirical risk minimization and portfolio selection problems, for which \cite{bertsimas2019sparse} and \cite{bertsimas2018portfolio} provide extensive numerical experiments.

\subsection{Overall Empirical Performance Versus State-of-the-Art}
In this section, we compare our approach to state-of-the-art methods, and demonstrate that our approach outperforms the state-of-the-art for several relevant problems.

\subsubsection{Network Design} \label{ssec:NDProb}
We begin by evaluating the performance of our approach for the multi-commodity network design problem \eqref{eqn:mcnd}. We adapt the methodology of \cite{gunluk2010perspective} and generate instances where each node $i \in [m]$ is the unique source of exactly one commodity ($k=m$). For each commodity $j \in [m]$, we generate demands according to
%\begin{align*}
{$b^{j}_{j'} = \lfloor \mathcal{U}(5, 25) \rceil$ for $j' \neq j$ %\mbox{ and }
and
$b^j_j = - \sum_{j' \neq j} b^{j}_{j'}$,}
%\end{align*}
where $\lfloor x \rceil$ is the closest integer to $x$ and $\mathcal{U}(a,b)$ is a uniform random variable on $[a,b]$. We generate edge construction costs, $c_e$, uniformly on $ \mathcal{U}(1, 4)$, and marginal flow circulation costs proportionally to each edge length\footnote{Nodes are uniformly distributed over the unit square $[0,1]^2$. We fix the cost to be ten times the Euclidean distance.}. The discrete set $\mathcal{Z}$ contains constraints of the form $\bm{z}_0 \leq \bm{z}$, where $\bm{z}_0$ is a binary vector which encodes existing edges. We generate graphs which contain a spanning tree plus $p m$ additional randomly picked edges, with $p \in [ 4 ]$, so that the initial network is {connected with $O(m)$ edges.} %, which is sparser and more realistic that $O(n^2)$ as considered by \cite{gunluk2010perspective}.
We also impose a cardinality constraint $\bm{e}^\top \bm{z} \leq (1+5\%) \bm{z}_0^\top \bm{e}$, which ensures that the network size increases by no more than $5\%$. For each edge, we impose a capacity $u_e \sim \lfloor \mathcal{U}(0.2, 1)B/A \rceil$, where $B = - \sum_{j=1}^m b_j^j$ is the total demand and $A = (1+p) m$. % is the number of existing edges.
We penalize the constraint $\bm{x} \leq \bm{u}$ with a penalty parameter $\lambda = 1,000$\footnote{We do so to allow for a fair comparison between big-$M$ and ridge regularization. By penalizing the capacity constraint, we remove a natural big-$M$ regularization term and no regularization can be considered as more natural than the other.}. For big-$M$ regularization, we set $M = \sum_j |b_j^j|$, and take $\gamma = \tfrac{2}{m(m-1)}$ for ridge regularization. % Note for Jean: can you add a sentence on the heuristics used for the warm-start for this set of experiments. Also can you please state the regularizer used to initially solve the problem.

We apply our approach to large networks with $100$s nodes, i.e., $10,000$s edges, which is ten times larger than the state-of-the-art \citep{holmberg1998solving,gunluk2010perspective}, and compare the quality of the incumbent solutions after an hour, {since no approach could terminate up to a satisfiable optimality gap within this time limit}.  Note that we define the quality of a solution as its cost in absence of regularization, although we might have augmented the original formulation with a regularization term to compute the solution. {As a result, we can compare the performance big-$M$ and ridge regularization directly, despite the fact that the optimization problems they solve are actually different.  On the other hand, performance metrics that depend on the function being minimized, such as the optimality gap, would not permit such a comparison.} In $100$ instances, our cutting plane algorithm with big-$M$ regularization provides a better solution $94\%$ of the time, by $9.9\%$ on average, and by up to $40\%$ for the largest networks.
For ridge regularization, the cutting plane algorithm scales to higher dimensions than plain mixed-integer SOCP, returns solutions systematically better than those found by CPLEX (in terms of unregularized cost), by $11\%$ on average. Also, ridge regularization usually outperforms big-$M$ regularization, as reported in Table \ref{tab:nd_big_objval}.
{Given how numerically challenging these optimization problems are, the optimality gaps returned by all methods are often uninformative ($>100 \%$) - see Section \ref{sec:A.nd} Table \ref{tab:nd_big_gap}. Still, we observe that, with big-$M$ regularization, CPLEX systematically returns tighter optimality gaps that the cutting-plane approach, while with ridge regularization, the gaps obtained by the cutting-plane algorithm are tighter $86\%$ of the times. All in all,}
even artificially added, ridge regularization improves the tractability of outer approximation. %the overall problem.

\begin{table} % Note for Jean: can you add a sentence on the lower bounds for this set of instances please
\centering
\footnotesize
\caption{Best solution found after one hour on network design instances with $m$ nodes and $(1+p)m$ initial edges. We report improvement, i.e., the relative difference between the solutions returned by CPLEX and the cutting-plane. Values are averaged over five randomly generated instances.
For ridge regularization, we report the ``unregularized'' objective value, that is we fix $\bm{z}$ to the best solution found and resolve the corresponding sub-problem with big-$M$ regularization. A ``$-$'' indicates that the solver could not finish the root node inspection within the time limit (one hour), and ``Imp.'' is an abbreviation of improvement.}
\begin{tabular}{llc|rrr|rrr|c}
& & & \multicolumn{3}{c}{Big-$M$} & \multicolumn{3}{c}{Ridge} & Overall\\
$m$ & $p$ & unit & CPLEX & Cuts & Imp. & CPLEX & Cuts & Imp. & Imp.\\
\toprule
40 & 0 & $ \times 10^{9} $ &  1.17 &  \textbf{1.16} & $ 0.86 \% $ & 1.55 &  \textbf{1.16} & $ 24.38 \% $ & $1.74\%$ \\
80 & 0 & $ \times 10^{9} $ &  8.13 &  7.52 & $ 6.99 \% $ & 9.95 &  \textbf{7.19} & $ 26.74 \% $  & $10.85\%$ \\
120 & 0 & $ \times 10^{10} $ &  3.03 &  2.10 & $ 29.94 \% $ & $-$ &  \textbf{1.94} & $ - \% $  & $35.30\%$ \\
160 & 0 & $ \times 10^{10} $ &  5.90 &  {4.32} & $ 26.69 \% $ & $-$ & \textbf{4.07} & $- \%$  & $30.91\%$ \\
200 & 0 & $ \times 10^{10} $ &  11.45 &  \textbf{7.78} & $ 31.45 \% $ & $-$ & \textbf{7.50} & $- \%$  & $32.32\%$  \\
% 240 & 0 & $ \times 10^{11} $ &  1.99 &  \textbf{1.19} & $ 39.22 \% $ & $-$ & $-$ & $- \%$  & $43.34\%$ \\
\midrule
40 & 1 & $ \times 10^{8} $ &  5.53 &  5.47 & $ 1.07 \% $ & 5.97 &  \textbf{5.45} & $ 8.74 \% $  & $1.41\%$ \\
80 & 1 & $ \times 10^{9} $ &  2.99 &  \textbf{2.94} & $ 1.81 \% $ & 3.16 &  2.95 & $ 6.78 \% $  & $1.89\%$ \\
120 & 1 & $ \times 10^{9} $ &  8.38 &  \textbf{7.82} & $ 6.69 \% $ & $-$ &  \textbf{7.82} & $-\% $  & $6.86\%$ \\
160 & 1 & $ \times 10^{10} $ &  1.64 &  \textbf{1.54} & $ 5.98 \% $ & $-$ &  \textbf{1.54} & $- \%$  & $6.03\%$ \\
200 & 1 & $ \times 10^{10} $ &  2.60 &  {2.54} & $ 2.33 \% $ & $-$ & \textbf{2.26} & $- \%$  & $12.98\%$ \\
% 240 & 1 & $ \times 10^{10} $ &  3.91 &  \textbf{2.79} & $ 28.81 \% $ & $-$ & $-$ & $- \%$  & $28.81\%$ \\
\midrule
40 & 2 & $ \times 10^{8} $ &  4.45 &  4.38 & $ 1.62 \% $ & 4.76 &  \textbf{4.36} & $ 8.27 \% $ & $2.06\%$  \\
80 & 2 & $ \times 10^{9} $ &  2.44 & \textbf{ 2.31} & $ 5.39 \% $ & 2.46 &  \textbf{2.31} & $ 5.97 \% $  & $5.40\%$ \\
120 & 2 & $ \times 10^{9} $ &  6.23 &  \textbf{5.89} & $ 5.55 \% $ & $-$ & \textbf{5.89} & $- \%$  & $5.75\%$ \\
160 & 2 & $ \times 10^{11} $ &  1.22 &  {1.16} & $ 4.74 \% $ & $-$ & \textbf{0.71} & $-\%$  & $19.33\%$ \\
200 & 2 & $ \times 10^{10} $ &  2.06 &  {1.43} & $ 30.46 \% $ & $-$ & \textbf{1.01} &  $- \%$  & $73.43\%$ \\
\midrule
40 & 3 & $ \times 10^{8} $ &  3.91 &  \textbf{3.85} & $ 1.58 \% $ & 4.13 &  \textbf{3.85} & $ 6.73 \% $  & $1.78\%$ \\
80 & 3 & $ \times 10^{9} $ &  2.06 & \textbf{1.94} & $ 5.76 \% $ & 2.04 &  \textbf{1.94} & $ 5.44 \% $  & $5.85\%$ \\
120 & 3 & $ \times 10^{9} $ &  5.43 &  {5.15} & $ 5.31 \% $ & $-$ & $\textbf{4.2}$ &  $- \%$  & $12.35\%$ \\ \midrule
40 & 4 & $ \times 10^{8} $ &  3.32 &  3.28 & $ 1.35 \% $ & 3.53 & \textbf{3.26} & $ 7.71 \% $  & $1.85\%$ \\
80 & 4 & $ \times 10^{9} $ &  1.88 &  \textbf{1.77} & $ 5.59 \% $ & $-$ & $ \textbf{1.77}$ & $- \%$  & $5.64\%$ \\
\bottomrule
\end{tabular}
\label{tab:nd_big_objval}
\end{table}

% \begin{table}
% \centering
% \footnotesize
% \caption{Optimality gap after one hour on network design instances with $m$ nodes and $(1+p)m$ initial edges. We only report results for instances where the resulting gap was less than $100\%$ for at least one of the four approaches. A ``$-$'' indicates that the solver could not finish the root node inspection within the time limit (one hour).   }
% \begin{tabular}{ll|rr|rr}
% & & \multicolumn{2}{c}{Big-$M$} & \multicolumn{2}{c}{Ridge} \\
% $m$ & $p$ & CPLEX & Cuts & CPLEX & Cuts \\
% \toprule
% 40 & 0 &  $69.8\%$ & $>100\%$  & $98.9\%$ & $96.7\%$ \\
% 80 & 0 &  $100.0\%$ &  $>100\%$ & $>100\%$ & $>100\%$ \\
% \midrule
% 40 & 1 &  $38.6\%$ &  $>100\%$ & $99.8\%$ & $97.1\%$ \\
% 80 & 1 &  $100.0\%$ &  $>100\%$ & $>100\%$ & $95.6\%$ \\
% 120 & 1 & $>100\%$ &  $>100\%$ & $-$ & $96.6\%$ \\
% \midrule
% 40 & 2 &  $23.3\%$ &  $>100\%$ & $>100\%$ & $97.7\%$ \\
% 80 & 2 &  $100.0\%$ &  $>100\%$ & $>100\%$ & $96.3\%$ \\
% \midrule
% 40 & 3 &  $74.6\%$ &  $>100\%$ & $97.5\%$ & $98.1\%$ \\
% \midrule
% 40 & 4 &  $100.0\%$ &  $>100\%$ & $99.2\%$ & $98.2\%$ \\
% 80 & 4 &  $100.0\%$ &  $>100\%$ & $>100\%$ & $80.2\%$ \\
% \bottomrule
% \end{tabular}
% \label{tab:nd_big_gap}
% \end{table}

\begin{table}
\centering
\caption{Average runtime in seconds on binary quadratic optimization problems from the Biq-Mac library \citep{wiegele2007biq, billionnet2007using}. Values are averaged over $10$ instances. A ``$-$'' denotes an instance which was not solved because the approach did not respect the $32$GB peak memory budget.}
\begin{tabular}{l l |r r r r@{}} \toprule
Instance & $n$  &  \multicolumn{4}{c}{Average runtime (s)/Average optimality gap $(\%)$} \\
\cmidrule{3-6} & &  CPLEX-M & CPLEX-M-Triangle &  Cuts-M & Cuts-M-Triangle \\ \midrule
bqp-$50$ & $50$ & $29.4$ & $0.6$ & $30.6$ & $\textbf{0.4}$ \\
bqp-$100$ & $100$ & $122.3$ & $51.7$ & $25.3\%$ & $\textbf{38.6}$ \\
bqp-$250$ & $250$ & $1108.1\%$ & $83.5\%$ & $87.0\%$ & $\textbf{46.1\%}$\\
bqp-$500$ & $500$ & $2055.8\%$ & $1783.3\%$ & $\textbf{157.3\%}$ & $410.7\%$\\
bqp-$1000$ & $1000$ & $-$  & $-$  & $\textbf{260.9}\%$ & $-$ \\
be$100$ & $100$ & $\textbf{79.7\%}$ & $208.0\%$ & $249.4\%$ & $201.2\%$ \\
be$120.8$ & $120$ & $\textbf{146.4\%}$ & $225.8\%$ & $264.1\%$ & $220.3\%$ \\
\bottomrule
\end{tabular}
\label{tab:bqp_aggregate_v1}
\end{table}

\subsubsection{Binary Quadratic Optimization} \label{ssec:bqpnumres}
We study some of the binary quadratic optimization problems collated in the BQP library by \cite{wiegele2007biq}. Specifically, the bqp-$\{50,100,250,500, 1000\}$ instances generated by \cite{beasley1990or}, which have a cost matrix density of $0.1$% and all non-zero entries drawn from $\mathcal{U}[-100, 100]$
, and the be-$100$ and be-$120.8$ instances generated by \cite{billionnet2007using}, which respectively have cost matrix densities of $1.0$ and $0.8$. %, non-zero on-diagonal entries drawn from $\mathcal{U}[-100, 100]$ and non-zero off-diagonal entries drawn from $\mathcal{U}[-50, 50]$.
Note that these instances were generated as maximization problems, and therefore we consider a higher objective value to be better.
We warm-start the cutting-plane approach with the best solution found after $10,000$ iterations of Goemans-Williamson rounding \citep[see][]{goemans1995improved}. We also consider imposing triangle inequalities \citep{deza2009geometry} via lazy callbacks, for they substantially tighten the continuous relaxations.

Within an hour, only the bqp-$50$ and bqp-$100$ instances could be solved by any approach considered here, in which case  cutting-planes with big-$M$ regularization is faster than CPLEX (see Table \ref{tab:bqp_aggregate_v1}). For instances which cannot be solved to optimality, although CPLEX has an edge in producing tighter optimality gaps for denser cost matrices, as depicted in Table \ref{tab:bqp_aggregate_v1}, the cutting-plane method provides tighter optimality gaps for sparser cost matrices, and provides higher-quality solutions than CPLEX for all instances, especially as $n$ increases (see Table \ref{tab:bqp_aggregate_ofv_1}).

We remark that the cutting plane approach has low peak memory usage compared with the other methods: For the bqp-$1000$ instances, cutting-planes without triangle inequalities was the only method which respected the $32$GB memory budget. This is another benefit of decomposing Problem \eqref{eqn:original_minlp} into master and sub-problems.
\begin{table}[h]
\centering
\caption{Average incumbent objective value (higher is better) after $1$ hour for medium-scale binary quadratic optimization problems from the Biq-Mac library \citep{wiegele2007biq, billionnet2007using}. ``$-$'' denotes an instance which was not solved because the approach did not respect the $32$GB peak memory budget. Values are averaged over $10$ instances. Cuts-Triangle includes an extended formulation in the master problem.}
\begin{tabular}{l l |r r r r@{}} \toprule
Instance & $n$  &  \multicolumn{4}{c}{Average objective value} \\
\cmidrule{3-6} & &  CPLEX-M & CPLEX-M-Triangle &  Cuts-M & Cuts-M-Triangle \\ \midrule
bqp-$250$ & $250$ & $9920.8$ & $41843.4$ & $\textbf{43774.9}$ & $43701.5$\\
bqp-$500$ & $500$ & $19417.1$ & $19659.0$ & $\textbf{122879.3}$ & $122642.4$\\
bqp-$1000$ & $1000$ & $-$  & $-$  & $\textbf{351450.7}$ & $-$ \\
be$100$ & $100$ & $16403.0$ & $16985.0$ & $17152.1$ & $\textbf{17178.5}$ \\
be$120.8$ & $120$ & $17943.2$ & $19270.3$ & $19307.7$ & $\textbf{19371.2}$ \\
\bottomrule
\end{tabular}
\label{tab:bqp_aggregate_ofv_1}
\end{table}

\subsubsection{Sparse Empirical Risk Minimization}  \label{ssec:serm}
\begin{figure}[h!]
\centering
\begin{subfigure}[t]{.47\linewidth}
	\centering
	\includegraphics[width=\linewidth]{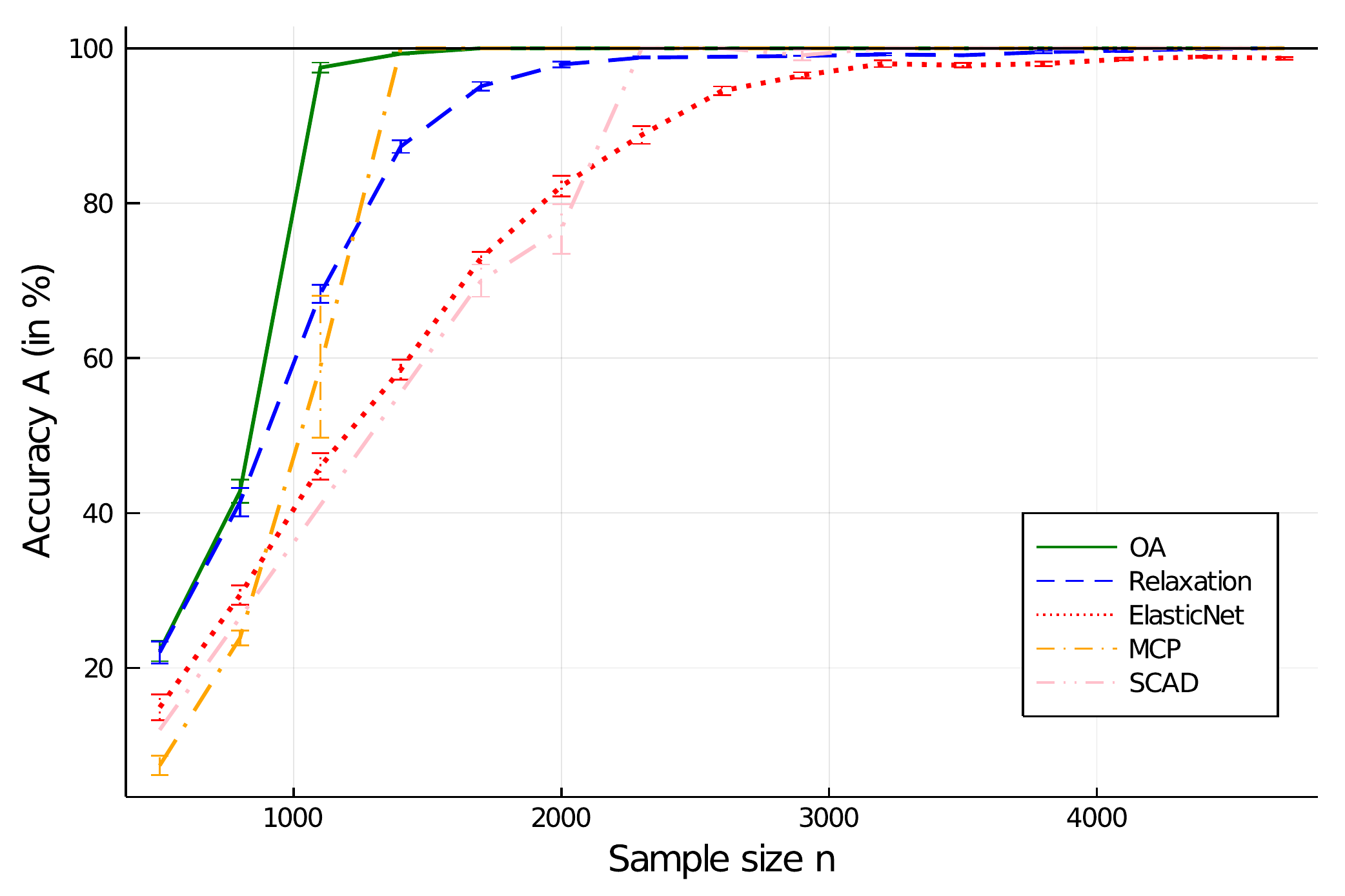}
	\caption{Regression, $p=20,000$}
\end{subfigure} %
\begin{subfigure}[t]{.47\linewidth}
	\centering
	\includegraphics[width=\linewidth]{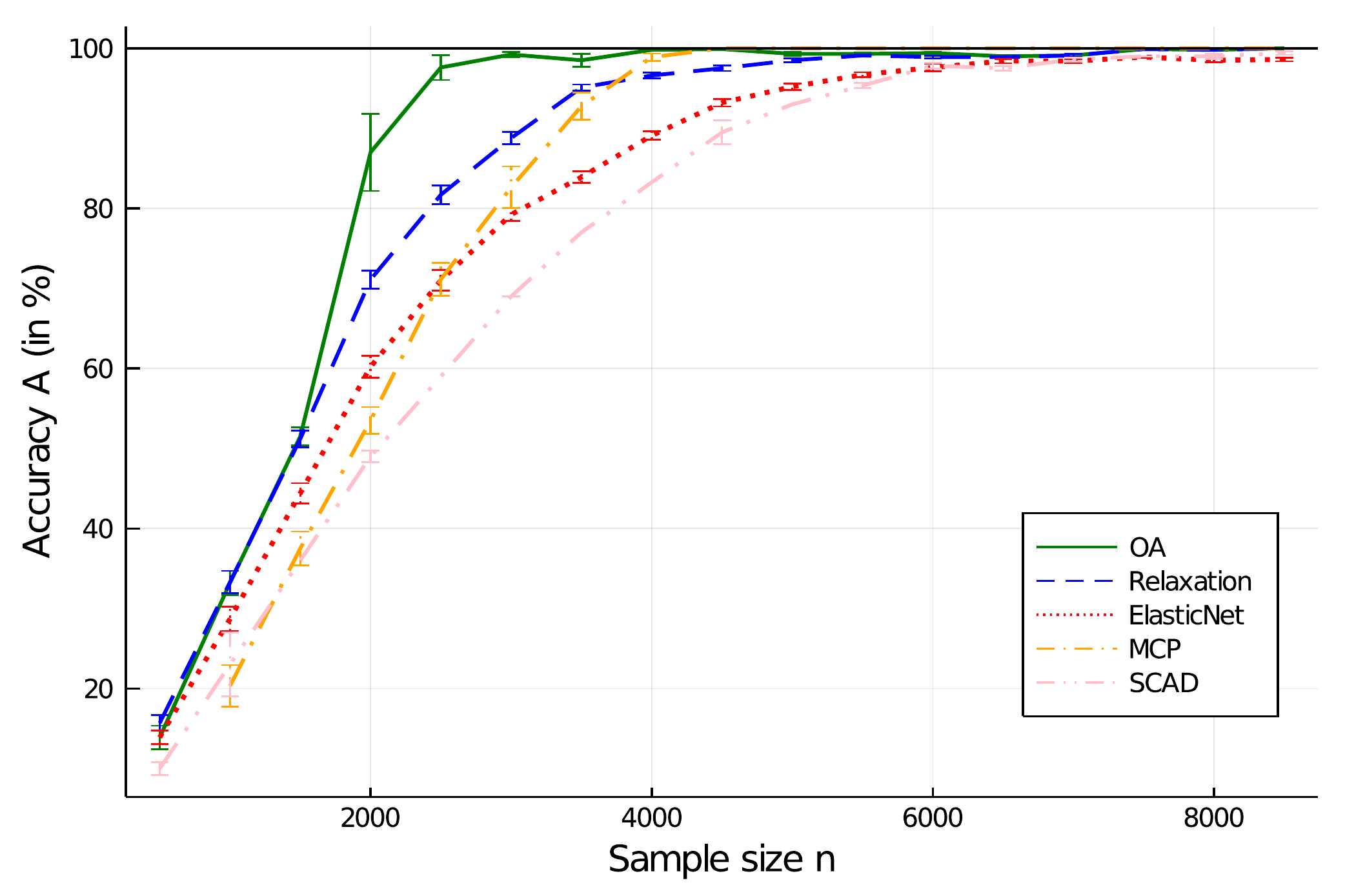}
	\caption{Classification, $p=10,000$}
\end{subfigure}
\caption{Accuracy ($A$) of the feature selection method as the number of samples $n$ increases, for the outer-approximation algorithm (in green), the solution found by the subgradient algorithm (in blue), ElasticNet (in red), MCP (in orange), SCAD (in pink) \citep[see][for definitions]{bertsimas2019sparse}. Results are averaged over $10$ instances of synthetic data with $(SNR, p, k) = (6, 20000, 100)$ for regression (left) and  $(5, 10000, 100)$ for classification (right).}
\label{fig:SERM_noMIC}
\end{figure}

For sparse empirical risk minimization, our method with ridge regularization scales to regression problems with up $p=100,000$s features and classification problems with $p=10,000$s of features \citep{bertsimas2019sparse}. This constitutes a three-order-of-magnitude improvement over previous attempts using big-$M$ regularization \citep{bertsimas2016best}. We also select features more accurately, as shown in Figure \ref{fig:SERM_noMIC}, which compares the accuracy of the features selected by the outer-approximation algorithm (in green) with those obtained from the Boolean relaxation (in blue) and other methods.
% \begin{table}[h]
% \centering\footnotesize
% \caption{Largest sparse regression instances reliably solved by each approach.}
% \begin{tabular}{@{}r l r @{}} \toprule
% Reference &  Solution method & Largest instance size solved \\
%  &  & (no. features) \\ \midrule
% %\cite{bertsimas2009algorithm} & Lemke pivot B$\&$B & $50$\\
% \cite{bertsimas2016best} & Big-$M$ B$\&$B & $100$s\\
% \cite{miyashiro2015mixed} & MISOCP B$\&$B & $100$\\
% \cite{pilanci2015sparse} & SOCP $+$ Randomized Rounding & (up to $15\%$ gap) $1000$s\\
% %\cite{dong2015regularization} & SDP $+$ Randomized Rounding & (up to $5\%$ gap) $100$s\\
% \textbf{\cite{bertsimas2017sparse}} & \textbf{Dual Branch-and-Cut} & $\bm{100,000}$\\
% %\cite{atamturk2019rank} & Rank$-1$ convexification & $100$\\
% \bottomrule
% \end{tabular}
% \label{tab:serm_summary}
% \end{table}
{
\subsubsection{Sparse Principal Component Analysis}
We applied our approach to sparse principal component analysis problems in \cite{bertsimas2020principal}, and by {(a) introducing either big-$M$ or ridge regularization and (b) introducing additional valid inequalities into the master problem, which we derived from the Gershgorin Circle Theorem \citep[see][Section 2.3, for details]{bertsimas2020principal}} successfully solved problems where $p=100$s to certifiable optimality, and problems where $p=1000$s to certifiable near optimality{, as reported in Table \ref{tab:spca}; we refer to \cite{bertsimas2020principal} for descriptions of the datasets studied and more extensive numerical experiments}. This constitutes an order-of-magnitude improvement over existing certifiably near-optimal approaches, which rely on semidefinite techniques and therefore cannot scale to $p=1000$s.

{
\begin{table}[h!]
\centering\footnotesize
\caption{Runtime in seconds per approach. We run all approaches on one thread, and impose a time limit of $600$s. If a solver fails to converge, we report the relative bound gap at termination in brackets, and the no. explored nodes and cuts at the time limit. For ridge regularization, we set $\gamma={100}/{k}$.}
\begin{tabular}{@{}l l l r r r r r r r@{}} \toprule
Dataset &  $p$ & $k$ & \multicolumn{3}{c@{\hspace{0mm}}}{Big-$M$ regularization} &  \multicolumn{3}{c@{\hspace{0mm}}}{Ridge regularization} \\
\cmidrule(l){4-6} \cmidrule(l){7-9} & &  & Time(s) & Nodes & Cuts & Time(s) & Nodes & Cuts \\\midrule
Pitprops & $13$ & $5$
& $\textbf{0.09}$ & $45$ & $22$ & $0.42$ & $42$ & $16$ \\
& & $10$
& $\textbf{0.08}$ & $223$ & $223$ & $0.68$ & $615$ & $244$  \\ \midrule
Wine & $13$ & $5$
& $\textbf{0.04}$ & $143$ & $69$ & $0.10$ & $73$ & $36$  \\\
& & $10$
& $\textbf{0.09}$ & $364$ & $232$ & $0.61$ & $394$ & $230$\\\midrule
Miniboone & $50$ & $5$
& $0.03$ & $3$ & $6$ & $\textbf{0.01}$ & $0$ & $2$ \\
&  & $10$
& $\textbf{0.04}$ & $4$ & $6$ & $0.07$ & $10$ & $13$ \\\midrule %alg 1 only off by 2% however
Communities & $101$ & $5$
& $\textbf{0.15}$ & $109$ & $2$ & $0.54$ & $272$ & $55$  \\
&  & $10$
& $\textbf{0.44}$ & $373$ & $76$ & $2.20$ & $1,800$ & $328$  \\\midrule
Arrhythmia & $274$ & $5$
& $\textbf{5.27}$ & $1,080$ & $192$ & $6.75$ & $1,242$ & $282$  \\
& & $10$
& ($\textbf{4.21\%}$) & $61,000$ & $11,600$ & ($4.63\%)$ & $77,200$ & $11,360$ \\\midrule
Micromass & $1300$ & $5$
& $\textbf{131.3}$ & $4,580$ & $4$ & $163.2$ & $4$ & $3,809$  \\
& & $10$
& $\textbf{378.6}$ & $321$ & $16,090$ & $510.3$ & $21,700$ & $566$ \\
\bottomrule
\end{tabular}
\label{tab:spca}
\end{table}

}
}

\subsubsection{Sparse Portfolio Selection} \label{ssec:sps}
We applied our approach to sparse portfolio selection problems in \cite{bertsimas2018portfolio}. By introducing a ridge regularization term, we successfully solved instances to optimality at a scale of one order of magnitude larger than previous attempts as summarized in Table \ref{tab:portfolio_summary}{. Specifically, we optimized over the securities in the Wilshire $5000$, which contains around $3,200$ securities, an improvement upon existing techniques, which cannot currently scale beyond the securities in the $S\&P$ 500. Moreover, at smaller scales which existing techniques have been benchmarked on—including the set of synthetic instances generated by \cite{frangioni2006perspective} with $200-400$ securities—our approach is as fast as and often faster than existing state-of-the-art approaches including \cite{zheng2014improving, frangioni2016approximated} among others \citep[see][Section 5.2, for details]{bertsimas2018portfolio}.}
\begin{table}[h!]
\centering\footnotesize
\caption{Largest sparse portfolio instances reliably solved by each approach}
\begin{tabular}{@{}r l r @{}} \toprule
Reference &  Solution method & Largest instance size solved\\
 &   & (no. securities) \\
\midrule
\cite{frangioni2009computational} & Perspective cut+SDP & $400$\\
\cite{bonami2009exact} & Nonlinear B$\&$B & $200$\\
\cite{gao2013optimal} & Lagrangian relaxation B$\&$B & $300$\\
\cite{cui2013convex} & Lagrangian relaxation B$\&$B & $300$\\
\cite{zheng2014improving} & SDP B$\&$B & $400$\\
\cite{frangioni2016approximated} & Approx. Proj. Perspective Cut & $400$\\
\cite{bertsimas2018portfolio} & Algorithm \ref{alg:OA} with ridge regularization & $3,200$\\
\bottomrule
\end{tabular}
\label{tab:portfolio_summary}
\end{table}

\subsection{Evaluation of Different Ingredients in Our Numerical Recipe}
%Individual contributions of our different ingredients}
\label{ssec:flpnumres}
We now consider the capacitated facility problem \eqref{eqn:flp} on 112 real-world instances available from the OR-Library \citep{beasley1990or,holmberg1999exact}, with the natural big-$M$ and the ridge regularization with $\gamma =1$. In both cases, the algorithms return the true optimal solution. Compared to CPLEX with big-$M$ regularization, our cutting plane algorithm with big-$M$ regularization is faster in $12.7\%$ of instances (by $53.6\%$ on average), and in $23.85\%$ of instances (by $54.5\%$ on average) when using a ridge penalty. This observation suggests that ridge regularization is better suited for outer-approximation, most likely because, as discussed in Section \ref{sec:algo.oa}, a strongly convex ridge regularizer breaks the degeneracy of the separation problems. Note that our approach could benefit from multi-threading and restarting.

We take advantage of these instances to breakdown the independent contribution of each ingredient in our numerical recipe in Table \ref{tab:flp_comp_breakdown}. Although each ingredient contributes independently, jointly improving the lower and upper bounds provides the greatest improvement.
\begin{table}[h!] \footnotesize
\centering
\caption{Proportion of wins and relative improvement over CPLEX in terms of computational time on the 112 instances from the OR-library \citep{beasley1990or,holmberg1999exact} for different implementations of our method: an outer-approximation (OA) scheme with cuts generated at the root node using Kelley's method (OA + Kelley), OA with the local search procedure (OA + Local search) and OA with a strategy for both the lower and upper bound (OA + Both). Relative improvement is averaged over all ``win'' instances.}
\begin{tabular}{l r r r r@{}} \toprule
%  &  &  \multicolumn{2}{c}{Runtime (s) / Average Bound Gap ($\%$)} \\
  & \multicolumn{2}{c}{Big-$M$} & \multicolumn{2}{c}{Ridge} \\
Algorithm & $\%$ wins & Relative improvement & $\%$ wins & Relative improvement \\
\midrule
OA + Kelley &$1.8\%$ & $36.6\%$ &$30.1\%$ & $91.6\%$\\
OA + Local search & $1.9\%$ & $49.5\%$ & $19.4\%$& $73.8\%$ \\
OA + Both & \textbf{$12.7\%$} & \textbf{$53.6\%$} & \textbf{$92.5\%$} & \textbf{$91.7\%$}\\\bottomrule
\end{tabular}
\label{tab:flp_comp_breakdown}
\end{table}

\subsection{Big-\texorpdfstring{$M$}{M} Versus Ridge Regularization} \label{ssec:UCnumres}
In this section, our primary interest is in ascertaining conditions under which it is advantageous to solve a problem using big-$M$ or ridge regularization, and argue that ridge regularization is preferable over big-$M$ regularization as soon as the objective is sufficiently strongly convex.

To illustrate this point, we consider large instances of the thermal unit commitment problem originally generated by \cite{frangioni2006solving}, and multiply the quadratic coefficient $a_i$ for each generator $i$ by a constant factor $\alpha \in \{0.1, 1, 2, 5, 10\}$. Table \ref{tab:uc_res_aggregate_cplex_v2} depicts the average runtime for CPLEX to solve both formulations to certifiable optimality, or provides the average bound-gap whenever CPLEX exceeds a time limit of $1$ hour. Observe that when $\alpha \leq 1$, the big-$M$ regularization is faster, but, when $ \alpha >1$ the MISOCP approach converges fast while the big-$M$ approach does not converge within an hour. Consequently, ridge regularization performs more favorably whenever the quadratic term is sufficiently strong.
\begin{table}[h!]
\centering\footnotesize
\caption{Average runtime in seconds per approach, on data from \cite{frangioni2006solving} where the quadratic cost are multiplied by a factor of $\alpha$. If the method did not terminate in one hour, we report the bound gap. $n$ denotes the number of generators, each instances has $24$ trade periods.}
\begin{tabular}{l| r r r r r r r r r r@{}} \toprule
%  &  &  \multicolumn{2}{c}{Runtime (s) / Average Bound Gap ($\%$)} \\
$\alpha$  & \multicolumn{2}{c}{$0.1$} & \multicolumn{2}{c}{$1$} & \multicolumn{2}{c}{$2$}& \multicolumn{2}{c}{$5$}& \multicolumn{2}{c}{$10$}\\
 \cmidrule{2-11}  $n$& Big-$M$ & Ridge & Big-$M$ & Ridge & Big-$M$ & Ridge & Big-$M$ & Ridge & Big-$M$ & Ridge\\\midrule
 $100$ & $\textbf{93.6}$ & $299.0$ & $\textbf{16.2}$ & $229.4$ & $0.32\%$ & $\textbf{47.9}$ & $1.68\%$ & $\textbf{4.6}$ & $2.76\%$ & $\textbf{6.0}$\\
 $150$ & $\textbf{35.6}$ & $352.1$ & $\textbf{6.2}$ & $28.3$ & $0.25\%$ & $\textbf{33.4}$ & $1.69\%$ & $\textbf{6.4}$ & $2.82\%$ & $\textbf{8.0}$\\
 $200$ & $\textbf{56.3}$ & $138.1$ & $\textbf{3.3}$ & $239.7$ & $0.24\%$ & $\textbf{112.9}$ & $1.62\%$ & $\textbf{16.7}$ & $2.81\%$ & $\textbf{21.2}$\\
\bottomrule
\end{tabular}
\label{tab:uc_res_aggregate_cplex_v2}
\end{table}

\begin{figure}[h]
\centering
\begin{subfigure}[t]{.47\linewidth}
	\centering
	\includegraphics[width=\linewidth]{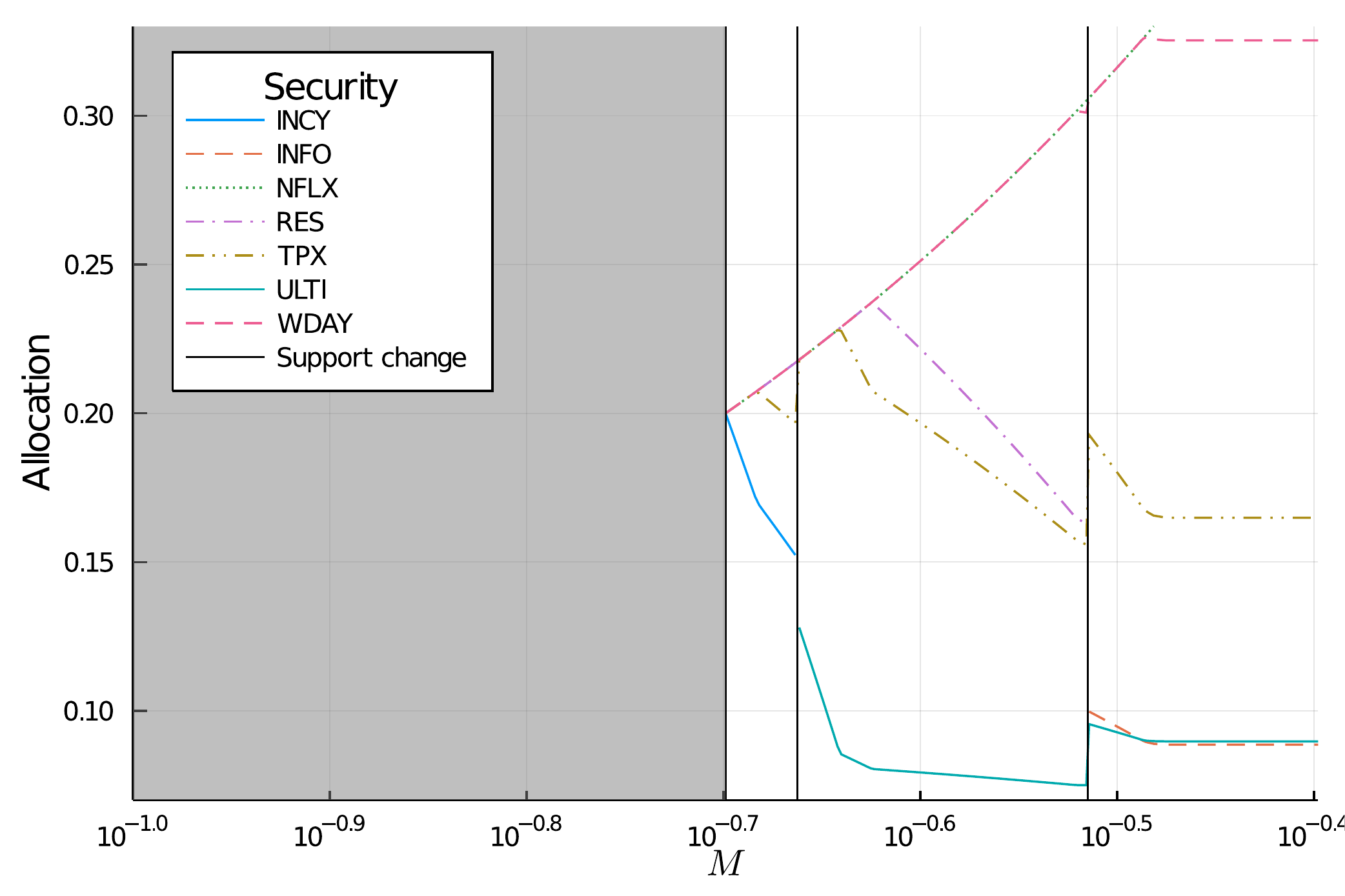}
	\caption{Big-$M$ regularization}
\end{subfigure} %
\begin{subfigure}[t]{.47\linewidth}
	\centering
	\includegraphics[width=\linewidth]{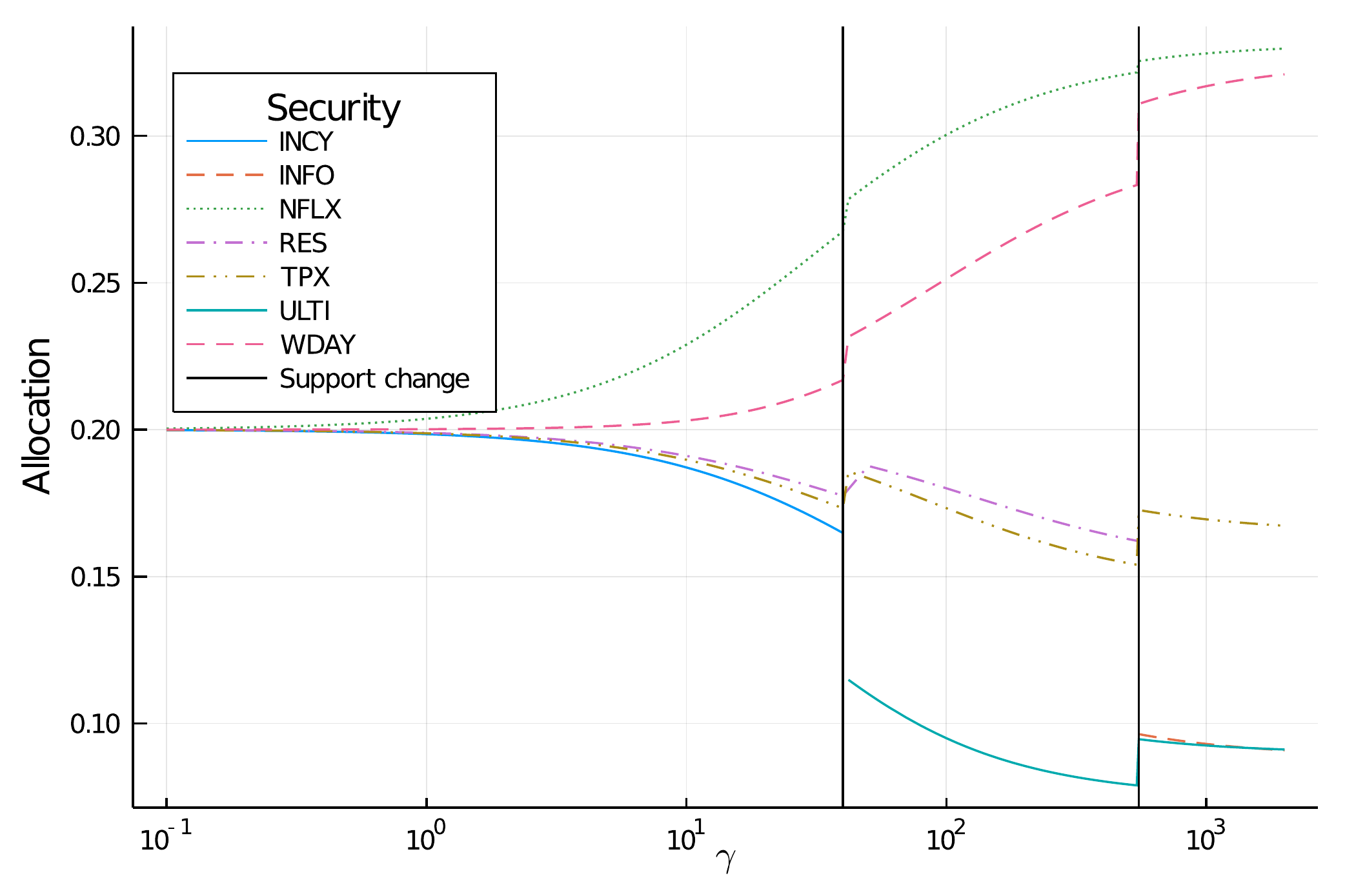}
	\caption{Ridge regularization}
\end{subfigure}
\caption{Optimal allocation of funds between securities as the regularization parameter ($M$ or $\gamma$) increases. Data is obtained from the Russell $1000$, with a cardinality budget of $5$, a rank$-200$ approximation of the covariance matrix, a one-month holding period and an Arrow-Pratt coefficient of $1$, as in \cite{bertsimas2018portfolio}. Setting $M <\frac{1}{k}$ renders the entire problem infeasible.}
\label{fig:sps_russ1000}
\end{figure}

\begin{figure}[h]
\centering
\begin{subfigure}[t]{.47\linewidth}
	\centering
	\includegraphics[width=\linewidth]{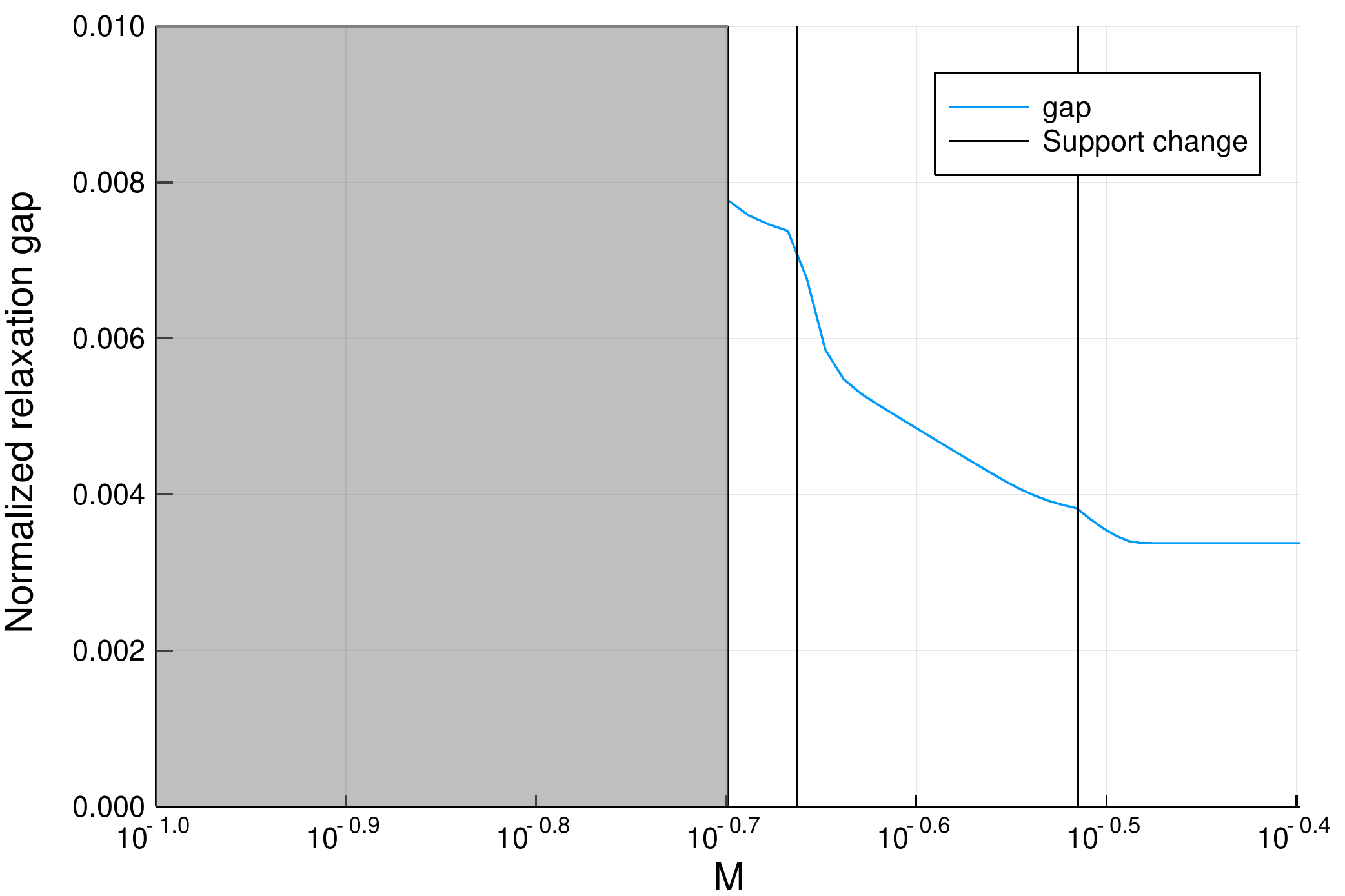}
	\caption{Big-$M$ regularization}
\end{subfigure} %
\begin{subfigure}[t]{.47\linewidth}
	\centering
	\includegraphics[width=\linewidth]{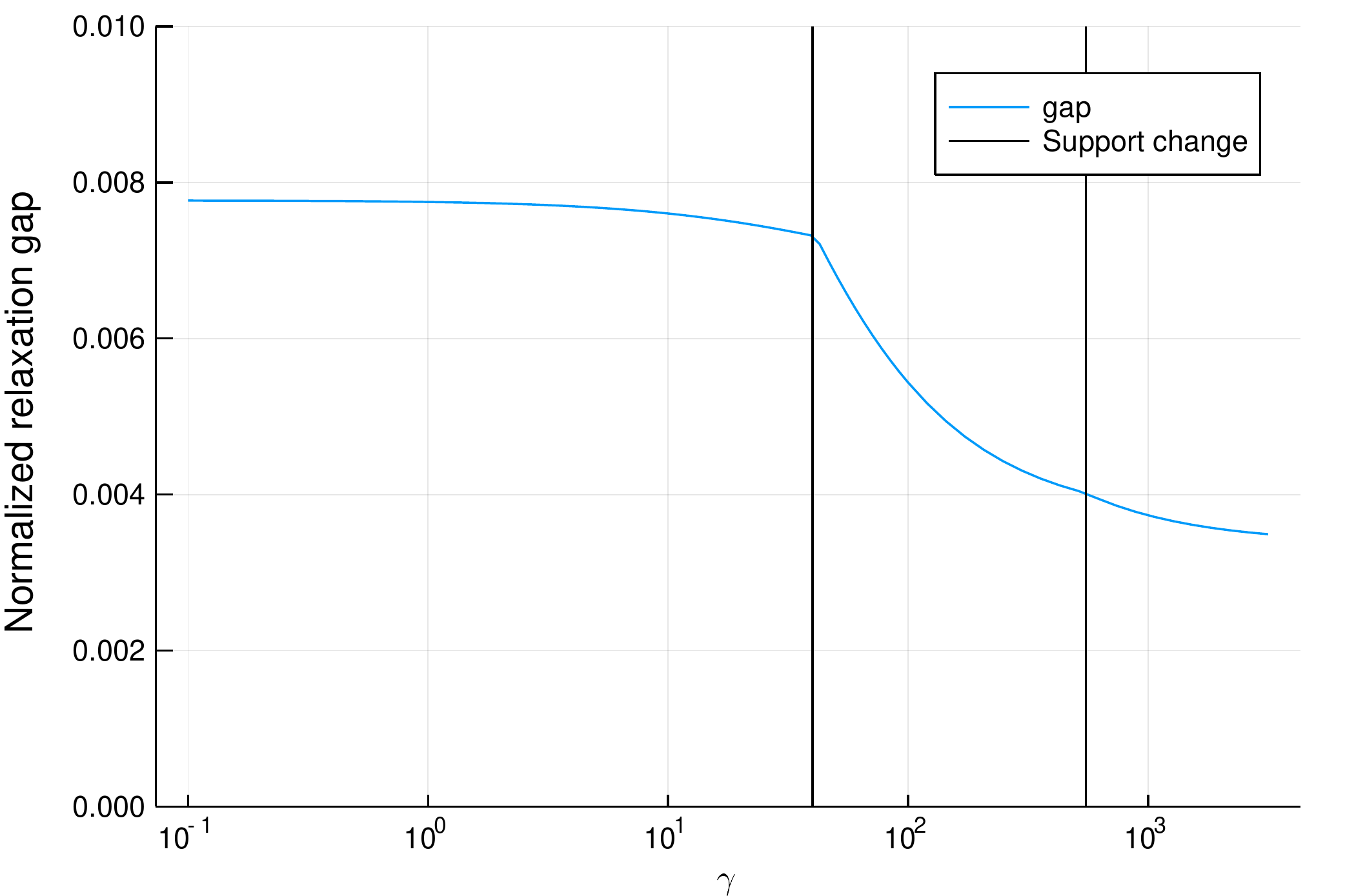}
	\caption{Ridge regularization}
\end{subfigure}
\caption{Magnitude of the normalized absolute bound gap as the regularization parameter ($M$ or $\gamma$) increases, for the portfolio selection problem studied in Figure \ref{fig:sps_russ1000}}
\label{fig:russ1000_gap}
\end{figure}

We also compare big-$M$ and ridge regularization for the sparse portfolio selection problem \eqref{eqn:sps}. Figure \ref{fig:sps_russ1000} depicts the relationship between the optimal allocation of funds $\bm{x}^\star$ and the regularization parameter $M$ (left) and $\gamma$ (right), and Figure \ref{fig:russ1000_gap} depicts the magnitude of the gap between the optimal objective and the Boolean relaxation's objective, normalized by the unregularized objective. The two investment profiles are comparable, selecting the same stocks. Yet, we observe two main differences: First, setting $M <\frac{1}{k}$ renders the entire problem infeasible, while the problem remains feasible for any $\gamma >0$. This is a serious practical concern in cases where a lower bound on the value of $M$ is not known a priori. Second, the profile for ridge regularization seems smoother than its equivalent with big-$M$.
% {
% Finally, we compare the size of the relative bound gaps (i.e., the gaps between the continuous relaxation and the optimal solution) under big-$M$ and ridge regularization for the network design problem with $m=10$, $p=0.2$ (averaged over $30$ instances). Figure \ref{fig:ndbgap} depicts the relationship between the size of the bound gap and the regularization parameter $M$ (left) and $\gamma$ (right). The magnitude of the bound gap remains stable at around $4.1\%$ for big-$M$ regularization, and varies between around $3.9\%$-$12\%$ for ridge regularization. However, setting $M<100$ renders at least one instance infeasible. }

\subsection{Relative Merits of Big-\texorpdfstring{$M$}{M}, Ridge Regularization: An Experimental Perspective}
{ We now conclude our comparison of big-$M$ and ridge regularization, as initiated in Sections \ref{ssec:meritstheory} and \ref{ssec:algmerits}, by indicating the benefits of big-$M$ and ridge regularization, from an experimental perspective:
\begin{itemize}
    \item As observed in Section \ref{ssec:UCnumres}, big-$M$ and ridge regularization play fundamentally the same role in reformulating logical constraints. This observation echoes our theoretical analysis in Section \ref{sec:framework}.
    \item As observed in the unit commitment and sparse portfolio selection problems studied in Section \ref{ssec:UCnumres}, ridge regularization should be the method of choice whenever the objective function contains a naturally occurring strongly convex term, which is sufficiently large.
    \item As observed for network design and capacitated facility location problems in sections \ref{ssec:NDProb}-\ref{ssec:flpnumres}, ridge regularization is usually more amenable to outer-approximation than big-$M$ regularization, because it eliminates most degeneracy issues associated with outer-approximating MINLOs.
    \item The efficiency of outer-approximation schemes relies on the speed at which separation problems are solved. In this regard, special problem-structure or cardinality constraints on the discrete variable $\bm{z}$ drastically help. This has been the case in network design, sparse empirical risk minimization and sparse portfolio selection problems in Section \ref{ssec:NDProb}.
\end{itemize}
}

\section{Conclusion}
In this paper, we proposed a new interpretation of the big-$M$ method, as a regularization term rather than a modeling trick. By expanding this regularization interpretation to include ridge regularization, we considered a wide family of relevant problems from the Operations Research literature and derived equivalent reformulations as mixed-integer saddle-point problems, which naturally give rise to theoretical analysis and computational algorithms. Our framework provides provably near-optimal solutions in polynomial time via solving Boolean relaxations and performing randomized rounding\footnote{By ``polynomial time'', we mean with respect to the dimensionality of the relaxation, assuming the relaxation is solved to a fixed and finite precision and can be described using the symmetric cones described by \cite{nesterov1994interior}, as occurs for all examples discussed in this paper.} as well as certifiably optimal solutions through an efficient branch-and-bound procedure, and indeed frequently outperforms the state-of-the-art in numerical experiments.

We believe our framework, which decomposes the problem into a discrete master problem and continuous subproblems, could be extended more generally to mixed-integer semidefinite optimization, as developed in \cite{bertsimas2019cosice, bertsimas2020principal}.
\FloatBarrier
\subsubsection*{Acknowledgments} We thank the associate editor and the two anonymous referees for their valuable comments which improved the paper.
{\footnotesize
\bibliographystyle{siamplain}

}

%\appendix
\FloatBarrier

\appendix

\section{Omitted Proofs}\label{sec:proof}
\subsection{Proof of Theorem \ref{thm:integrality.gap}: Quality of the Random Rounding Strategy} \label{sec:proof.randomrounding}

\proof
We only detail the proof for the big-$M$ regularization case, as the ridge regularization case follows \textit{mutatis mutandis}. From Proposition \ref{thm:lipschitz},
\begin{align*}
0 \leq f(\bm{z})-f(\bm{z}^\star) \leq M  L |\mathcal{R}| \max_{\bm{\alpha} \geq \bm{0} : \| \bm{\alpha} \|_1 \leq 1} \sum_{i\in \mathcal{R}} (z^\star_i - z_i) \alpha_i.
\end{align*}
The polyhedron $\{ \bm{\alpha} : \bm{\alpha} \geq \bm{0}, \| \bm{\alpha} \|_1 \leq 1\}$ admits $|\mathcal{R}|+1$ extreme points. However, if
\begin{align*}
\max_{\bm{\alpha} \geq \bm{0} : \| \bm{\alpha} \|_1 \leq 1} \sum_{i \in \mathcal{R}} (z^\star_i - z_i) \alpha_i > t,
\end{align*}for some $t>0$, then the maximum can only occur at some $\alpha > \mathbf{0}$ so that we can restrict our attention to the $|\mathcal{R}|$ positive extreme points. Applying tail bounds on the maximum of sub-Gaussian random variables over a polytope \citep[see][Theorem 1.16]{rigollet2015high}, since $\| \bm{\alpha} \|_2 \leq \| \bm{\alpha} \|_1 \leq 1$, we have for any $t>0$,
\begin{align*}
\mathbb{P}\left( \max_{\bm{\alpha} \geq \bm{0} : \| \bm{\alpha} \|_1 \leq 1} \sum_{i \in \mathcal{R}} (z^\star_i - z_i) \alpha_i > t \right) \leq | \mathcal{R} | \exp \left( - \dfrac{t^2}{2} \right),
\end{align*}
so that
\begin{align*}
\mathbb{P}\left( M L |\mathcal{R}| \max_{\bm{\alpha} \geq \bm{0} : \| \bm{\alpha} \|_1 \leq 1} \sum_{i \in \mathcal{R}} (z^\star_i - z_i) \alpha_i > \varepsilon \right) \leq | \mathcal{R} | \exp \left( - \dfrac{\varepsilon^2}{2 M^2 L^2 |\mathcal{R}|^2} \right).
\end{align*}
\endproof

\subsection{Proof of Theorem \ref{thm:perspective}: Relationship With Perspective Cuts}\label{ssec:proofofperspectivethm}

\proof Let us fix $\bm{z} \in \mathcal{Z}$. Then, we have that:
\begin{align*}
\max_{\bm{\alpha}} \: h(\bm{\alpha} ) - \dfrac{\gamma}{2} \sum_{j=1}^n z_j \alpha_j^2
&= \max_{\bm{\alpha}, \bm{\beta}} \: h(\bm{\alpha}) - \dfrac{\gamma}{2} \sum_{j=1}^n z_j \, {\beta}_j^2
\mbox{  s.t.  } \bm{\beta} = \bm{\alpha}, \\
&= \max_{\bm{\alpha}, \bm{\beta}} \: \min_{\bm{x}} \: h(\bm{\alpha}) - \dfrac{\gamma}{2} \sum_{j=1}^n z_j \, {\beta}_j^2 - \bm{x}^\top (\bm{\beta} - \bm{\alpha}), \\
&= \min_{\bm{x}} \: \underbrace{ \max_{\bm{\alpha}} \:  \left[ h(\bm{\alpha}) + \bm{x}^\top \bm{\alpha} \right]}_{(-h)^\star(\bm{x}) = g(\bm{x})} + \sum_{i=1}^n \max_{\beta_i} \: \left[ - \dfrac{\gamma}{2} z_i \, \beta_i^2 - x_i \beta_i \right].
\end{align*}
Finally, observing that
\begin{align*}
\max_{\beta_i} \: \left[ - \dfrac{\gamma}{2} z_i \, \beta_i^2 - x_i \beta_i \right] &= \begin{dcases}
\dfrac{x_i^2}{2\gamma z_i} &\mbox{  if } z_i > 0, \\
\max_{\beta_i} \: x_i \beta_i &\mbox{  if  } z_j = 0,
\end{dcases}
\end{align*}
concludes the proof.
\endproof

\section{Bounding the Lipschitz Constant}\label{sec:lipcont} %of $f(\bm{z})$
In our results, we relied on the observation that there exists some constant $L>0$ such that, for any $\bm{z} \in \mathcal{Z}$, $\| \bm{\alpha}^\star(\bm{z}) \| \leq L$. Such an $L$ always exists, since $\mathcal{Z}$ is a finite set. However, as our randomized rounding results depend on $L$, explicit bounds on $L$ are desirable.

We remark that while our interest is in the Lipschitz constant with respect to ``$\bm{\alpha}$'' in a generic setting, we have used different notation for some of the problems which fit in our framework, in order to remain consistent with the literature. In this sense, we are also interested in obtaining a Lipschitz constant with respect to $\bm{w}$ for the portfolio selection problem \eqref{eqn:sps}, among others.

In this appendix, we bound the magnitude of $L$ in a less conservative manner. Our first result provides a bound on $L$ which holds whenever the function $h(\bm{\alpha})$ in Equation \eqref{eqn:saddlepointproblem} is strongly concave in $\bm{\alpha}$, which occurs for the sparse ERM problem \eqref{eqn:serm} with ordinary least-squares loss, the unit commitment problem \eqref{eqn:UC}, the portfolio selection \eqref{eqn:sps}, and network design problems whenever $\bm{\Sigma}$ (resp. $\bm{Q}$) is full-rank:
\begin{lemma}\label{lemma:lconstant}
Let $h(\cdot)$ be a strongly concave function with parameter $\mu>0$ \citep[see][Chapter 9.1.2 for a general theory of strong convexity]{boyd2004convex}, and suppose that $\bm{0} \in \text{dom}(g)$ and $\bm{\alpha}^\star:=\arg \max_{\bm{\alpha}} \, h(\bm{\alpha})$. Then, for any choice of $\bm{z}$, we have
\begin{align*}
    \| \bm{\alpha}^\star(\bm{z}) \|_2^2 \leq 8 \, \frac{h(\bm{\alpha}^\star)-h(\bm{0})}{\mu},
\end{align*}
i.e., $\Vert \bm{\alpha}^\star(\bm{z})\Vert_\infty \leq L$, where $L:=2\sqrt{2\frac{h(\bm{\alpha}^\star)-h(\bm{0})}{\mu}}$.
\end{lemma}
\proof
By the definition of strong concavity, for any $\bm{\alpha}$ we have $$h(\bm{\alpha}) \leq h(\bm{\alpha}^\star) + \nabla h(\bm{\alpha}^\star)^\top (\bm{\alpha}-\bm{\alpha}^\star)  -\frac{\mu}{2}\Vert \bm{\alpha}-\bm{\alpha}^\star \Vert_2^2,$$ where $\nabla h(\bm{\alpha}^\star)^\top (\bm{\alpha}-\bm{\alpha}^\star) \leq {0}$ by the first-order necessary conditions for optimality, leading to
$$\Vert \bm{\alpha} - \bm{\alpha}^\star \Vert_2^2 \leq 2 \, \frac{h(\bm{\alpha}^\star)-h(\bm{\alpha})}{\mu}.$$
In particular for $\bm{\alpha} = \bm{0}$, we have
$$\Vert \bm{\alpha}^\star \Vert_2^2 \leq 2 \, \frac{h(\bm{\alpha}^\star)-h(\bm{0})}{\mu},$$
and for $\bm{\alpha} = \bm{\alpha}^\star(\bm{z})$,
$$\Vert \bm{\alpha}^\star(\bm{z})-\bm{\alpha}^\star \Vert_2^2 \leq 2 \,\frac{h(\bm{\alpha}^\star)-h(\bm{0})}{\mu},$$
since
$$ h(\bm{\alpha}^\star({\bm{z}})) \geq h(\bm{\alpha}^\star(\bm{z}))- \sum_{j=1}^n z_j \Omega_j^\star({\alpha}^\star (\bm{z})_j) \geq h(\bm{0}).$$ The result then follows by the triangle inequality.
\endproof

An important special case of the above result arises for the sparse ERM problem, as we demonstrate in the following corollary to Lemma \ref{lemma:lconstant}:
\begin{corollary}
For the sparse ERM problem \eqref{eqn:serm} with an ordinary least squares loss function and a cardinality constraint $\bm{e}^\top \bm{z} \leq k$, a valid bound on the Lipschitz constant is given by \begin{align*}\Vert \bm{\beta}^\star(\bm{z}) \Vert_\infty=\Vert \mathrm{Diag}(\bm{Z})\bm{X}^\top \bm{\alpha}^\star(\bm{z})\Vert_\infty & \leq \Vert \mathrm{Diag}(\bm{Z})\bm{X}^\top\Vert_\infty \Vert \bm{\alpha}^\star(\bm{z})\Vert_\infty \\ & \leq \max_{i}\bm{X}_{i, [k]} \Vert \bm{\alpha}\Vert_2 \leq 2\max_{i}\bm{X}_{i, [k]} \Vert \bm{y}\Vert_2,\end{align*} where $\bm{X}_{i, [k]}$ is the sum of the $k$ largest entries in the column $\bm{X}_{i, [k]}$.
\end{corollary}
\proof
Applying Lemma \ref{lemma:lconstant} yields the bound
$$ \Vert \bm{\alpha}\Vert_2 \leq 2 \Vert \bm{y}\Vert_2,$$
after observing that we can parameterize this problem in terms of $\bm{\alpha}$, and for this problem:
\begin{enumerate}
    \item Setting $\bm{\alpha}=0$ yields $h(\bm{\alpha})=0$.
    \item $0 \leq h(\bm{\alpha}^\star) \leq \bm{y}^\top \bm{\alpha}^\star-\frac{1}{2}\bm{\alpha}^{\star\top} \bm{\alpha}^\star \leq \frac{1}{2}\bm{y}^\top \bm{y}$.
    \item $h(\cdot)$ is strongly concave in $\bm{\alpha}$, with concavity constant $\mu\geq 1$.
\end{enumerate}
The result follows by applying the definition of the operator norm, and pessimizing over $\bm{z}$.
\endproof
% \subsection{A Bound from the Fenchel-Young Inequality}
% \begin{lemma}
% Let $e_i \in \mathrm{dom}(g)$, and let $h(\bm{\alpha})$ be non-negative on its domain. Then, for Equation \eqref{eqn:saddlepointproblem}, a valid Lipschitz bound on $\bm{\alpha}^\star(\bm{z})$ is given by $$\vert \alpha_i \vert \leq \max\left(g(\bm{e}_i), g(-\bm{e}_i) \right) $$
% \end{lemma}
% \begin{proof}
% By the definition of $h(\bm{\alpha})$, for any feasible $\bm{v}$ we have that:
% \begin{align*}
% h(\bm{\alpha}) \leq g(\bm{v})-\bm{v}^\top \bm{\alpha}.
% \end{align*}
% Re-arranging then reveals that:
% \begin{align*}
%     \bm{\alpha}^\top \bm{v} \leq g(\bm{v})-h(\bm{\alpha}).
% \end{align*}
%     Finally, substituting $\bm{e}_i=\bm{v}$ yields the result, as $h(\bm{\alpha}) \geq 0$.
% \end{proof}
{
\section{Supplementary material for the numerical experiments} \label{sec:A.nd}
In this section, we report additional performance metrics for the network design experiments presented in Section \ref{ssec:NDProb}. There, we reported the quality of the solution returned by all methods within one hour, and compared two regularization strategies (big-$M$ vs. ridge) and two algorithms (CPLEX vs. Cuts). Indeed, given the sizes of problems considered, the network design instances are computationally very challenging to solve. At such scales, finding a good feasible solution is already a very difficult task. In practice, this translates into optimality gaps that are often irrelevant (i.e., higher than $100\%$) in most of the instances. Table \ref{tab:nd_big_gap} reports the optimality gaps returned by each method after one hour, on instances where at least one of the four gaps was less than $100\%$.

\begin{table}
\centering
\footnotesize
\caption{Optimality gap after one hour on network design instances with $m$ nodes and $(1+p)m$ initial edges. We only report results for instances where the resulting gap was less than $100\%$ for at least one of the four approaches. A ``$-$'' indicates that the solver could not finish the root node inspection within the time limit (one hour).   }
\begin{tabular}{ll|rr|rr}
& & \multicolumn{2}{c}{Big-$M$} & \multicolumn{2}{c}{Ridge} \\
$m$ & $p$ & CPLEX & Cuts & CPLEX & Cuts \\
\toprule
40 & 0 &  $69.8\%$ & $100\%$  & $98.9\%$ & $96.7\%$ \\
80 & 0 &  $100\%$ &  $100\%$ & $100\%$ & $100\%$ \\
\midrule
40 & 1 &  $38.6\%$ &  $100\%$ & $99.8\%$ & $97.1\%$ \\
80 & 1 &  $100\%$ &  $100\%$ & $100\%$ & $95.6\%$ \\
120 & 1 & $100\%$ &  $100\%$ & $-$ & $96.6\%$ \\
\midrule
40 & 2 &  $23.3\%$ &  $100\%$ & $>100\%$ & $97.7\%$ \\
80 & 2 &  $100\%$ &  $100\%$ & $100\%$ & $96.3\%$ \\
\midrule
40 & 3 &  $74.6\%$ &  $>100\%$ & $97.5\%$ & $98.1\%$ \\
\midrule
40 & 4 &  $100\%$ &  $100\%$ & $99.2\%$ & $98.2\%$ \\
80 & 4 &  $100\%$ &  $100\%$ & $100\%$ & $80.2\%$ \\
\bottomrule
\end{tabular}
\label{tab:nd_big_gap}
\end{table}
}

\end{document}